\def\versionno{ Z2 = BFRS2  --  version  7.4   -- by tb --   2.2.10   }

\documentclass[12pt]{article}
\usepackage{amsmath, amsthm, amsfonts, enumerate, amssymb, bbm, xspace}
\usepackage[arrow, matrix, curve, all, knot]{xy}
\usepackage{graphics}
\usepackage{rotating}
\usepackage[english]{babel}

\usepackage{ifpdf} 
\ifpdf
\usepackage{epstopdf}
\usepackage{hyperref}
\else
\usepackage[hypertex]{hyperref}
\fi

\newtheorem{theorem}{Theorem}
\newtheorem{lemma}[theorem]{Lemma}
\newtheorem{conv}[theorem]{Convention}
\newtheorem{proposition}[theorem]{Proposition}
\newtheorem{corollary}[theorem]{Corollary}
\theoremstyle{definition}
\newtheorem{rem}[theorem]{Remark}
\newtheorem{definition}[theorem]{Definition}

\DeclareMathOperator{\End}{End}
\DeclareMathOperator{\Hom}{Hom}

 \numberwithin{equation}{section}

\catcode`\@=11
\newif\if@fewtab\@fewtabtrue
{\count255=\time\divide\count255 by 60
\xdef\hourmin{\number\count255}
\multiply\count255 by-60\advance\count255 by\time
\xdef\hourmin{\hourmin:\ifnum\count255<10 0\fi\the\count255}}
\def\ps@draft{\let\@mkboth\@gobbletwo
    \def\@oddfoot{\hbox to 7 cm{\tiny \versionno
       \hfil}\hskip -7cm\hfil\rm\thepage \hfil {\tiny\draftdate}}
    \def\@oddhead{}
    \def\@evenhead{}\let\@evenfoot\@oddfoot}
\def\draftdate{\number\month/\number\day/\number\year\ \ \ \hourmin }

\catcode`\@=12

\def\be            {\begin{equation}}
\def\bearl         {\begin{array}{l}} 
\def\bearll        {\begin{array}{ll}}
\def\bearlll       {\begin{array}{lll}}
\def\boti          {\,{\boxtimes}\,}
\def\botih         {\,{\widehat\boxtimes}\,}

\def\C             {\ensuremath{\mathcal C}\xspace}

\def\calI          {\mathcal{I}}
\def\calC          {\C}
\def\CbC           {\ensuremath{\hspace{.3pt}\mathcal C\hspace{.6pt}{\boxtimes}%
                   \hspace{1.4pt}\mathcal C}\xspace}
\def\CC            {\ensuremath{\mathcal C{\boxtimes}\mathcal C}\xspace}
\def\cir           {\,{\circ}\,}
\def\complex       {{\mathbb C}}
\def\congTo        {\,{\stackrel\cong\to}\,} 
\def\CtC           {\ensuremath{\hspace{.3pt}\mathcal C\hspace{.4pt}{\times}%
                   \hspace{1.1pt}\mathcal C}\xspace}
\def\dim           {\mathrm{dim}}
\def\Dim           {\mathrm{Dim}}
\def\dimc          {\mathrm{dim}_\complex}

\def\dsty          {\displaystyle }
\def\ee            {\end{equation}}
\def\eear          {\end{array}}
\def\unit          {\mathbf{1}}  
\newcommand\Endfun[2]{\mathcal{E}\hspace{-1.6pt}{\it nd}_{#1}(#2)}  
\def\eps           {\varepsilon}
\def\eq            {\,{=}\,}
\newcommand\erf[1] {\mbox{(\ref{#1})}}
\def\ev            {ev_{\unit,\unit}}
\def\gamm          {\gamma}
\def\I             {\calI} 
\def\id            {\mbox{\sl id}}

\def\iN            {\,{\in}\,} 
\def\intEnd        {\underline{\End}}
\def\mortimes	   {{\,\otimes_\complex\,}}
\def\M             {\ensuremath{\mathcal M}\xspace}
\def\N             {\ensuremath{\mathcal N}\xspace}
\newcommand\NN[3]  {{N_{#1#2}}^{\!\!#3}}
\def\obj           {{\mathcal O}bj}
\def\one           {\unit}

\def\oti           {\,{\otimes}\,}
\def\Oti           {{\otimes}}
\def\To            {\,{\to}\,}

\def\vectk         {{\mathcal V}\mbox{\sl ect}_\Bbbk}
\def\X             {\ensuremath{\mathcal X}}
\def\XBB           {\ensuremath{\X_{\!B|B}}}
\def\V             {\ensuremath{\mathcal V}}
\def\Vee           {{}^{\!\vee}}
\def\zet           {\mathbb Z}

\newcommand\bild[1]{\includegraphics{#1}}
\def\bg         {\begin{gather}}
\def\eg         {\end{gather}}
\def\bp         {\begin{picture}}
\def\ep         {\end{picture}}
\renewcommand\epsilon{\varepsilon}

\newcommand\void[1]{}

\begin{document}

\thispagestyle{empty}
\def\thefootnote{\fnsymbol{footnote}}
\begin{flushright}
   {\sf KCL-MTH-08-08}\\
   {\sf ZMP-HH/08-17}\\
   {\sf Hamburger$\;$Beitr\"age$\;$zur$\;$Mathematik$\;$Nr.$\;$320}\\[2mm]
   Dezember 2008
\end{flushright}
\vskip 2.0em
\begin{center}\Large 
  MODULE CATEGORIES \\[2pt] FOR PERMUTATION MODULAR INVARIANTS

\end{center}\vskip 1.4em
\begin{center}
  Till Barmeier\,$^{\,a}$,~
  ~J\"urgen Fuchs\,$^{\,b}$,~
  ~Ingo Runkel\,$^{\,c}$,~
  ~Christoph Schweigert\,$^{\,a}$\footnote{\scriptsize 
  ~Email addresses: \\
  $~$\hspace*{2.4em}barmeier@math.uni-hamburg.de, jfuchs@fuchs.tekn.kau.se, 
  ingo.runkel@kcl.ac.uk, schweigert@math.uni-hamburg.de}
\end{center}

\vskip 3mm

\begin{center}\it$^a$
  Organisationseinheit Mathematik, \ Universit\"at Hamburg\\
  Schwerpunkt Algebra und Zahlentheorie\\
  Bundesstra\ss e 55, \ D\,--\,20\,146\, Hamburg
\end{center}
\begin{center}\it$^b$
  Teoretisk fysik, \ Karlstads Universitet\\
  Universitetsgatan 21, \ S\,--\,651\,88\, Karlstad
\end{center}
\begin{center}\it$^c$
  Department of Mathematics, King's College London \\
  Strand, UK\,--\,London WC2R 2LS
\end{center}
\vskip 3.5em

\begin{abstract}
We show that a braided monoidal category \C\ can be endowed with the structure 
of a right (and left) module category over $\C \,{\times}\, \C$. In fact, 
there is a family of such module category structures, and they are mutually 
isomorphic if and only if \C\ allows for a twist. For the case that
\C\ is premodular we compute the internal End of the tensor unit of \C, 
and we show that it is an Azumaya algebra if \C\ is modular. As an application to
two-dimensional rational conformal field theory, we show that the module 
categories describe the permutation modular invariant for models based
on the product of two identical chiral algebras. It follows in particular that
all permutation modular invariants are physical.
\end{abstract}

\newpage

\section{Introduction}

Just as a monoidal category \X\ can be understood as the categorification of 
a ring, a right module category \M\ is the categorification of a right module 
over a ring. Namely, there is a bifunctor $\boxtimes{:}\ \M\,{\times}\,\X \To \M$ 
together with associativity and unit constraints obeying suitable compatibility 
conditions (they will be recalled in section \ref{sec:familiy} below). 
Module categories are for instance used to define Morita equivalence between 
monoidal categories \cite{muge8}, and they play a significant role in the study 
of two-dimensional rational conformal field theory \cite{scfr2}.

\medskip

Any ring can be regarded as a right module over itself by setting 
$m\,{.}\,a\,{:=}\,m\,{\cdot}\, a$. In the same way, the tensor product and the 
associativity constraints can be used to endow any monoidal category \C\ with the 
structure of a right module category over itself. A {\em commutative\/} ring $R$ 
carries the structure of a right module over the ring $R\,{\otimes_\zet}\, R$
by $a\,{.}\,(b\oti c)\,{:=}\, a\,{\cdot}\,b\,{\cdot}\,c$. Notice that the module property
requires the equality
\begin{equation*}
\begin{split}
  m \,{\cdot}\, a_1 \,{\cdot}\, a_2 \,{\cdot}\, b_1 \,{\cdot}\, b_2 
  &\equiv m \,{.} \left((a_1\oti b_1)\,{\cdot}\, (a_2\oti b_2)\right)\\
  &= \left(m \,{.}\, (a_1\oti b_1)\right) {.}\, (a_2\oti b_2)
  \equiv m\,{\cdot}\, a_1\,{\cdot}\, b_1\,{\cdot}\, a_2\,{\cdot}\, b_2 \,,
\end{split}
\end{equation*}
which holds by commutativity of $R$. The categorification of a commutative 
ring is an additive braided monoidal category. Hence the question arises whether any
braided monoidal category \C\ carries the structure of a right module category over 
the monoidal category \CtC.  We show that this is indeed the case. However,
the lift of the permutation to an element 
in the braid group is not unique. Using inequivalent lifts, we present a family
of module category structures, indexed by integers, on these categories and 
show that they are mutually equivalent if and only if \C\ can be equipped with 
a twist, or balancing, see Theorems \ref{thm_onlybraided} and \ref{thm_psi-twist}
and Remark \ref{rem_G_gamma-vs-id_gamma}.

Given a module category \M\ over \X, denote by $\Endfun\X\M$ the category
of module endofunctors of \M. This is a monoidal category, with composition 
as tensor product. We recall in section \ref{sec:braid} that if \X\ is 
braided there are two monoidal functors $\alpha^\pm{:}~ \X \To \Endfun\X\M$ 
called braided induction. In section \ref{sec:azumaya} we show that
if \C\ is a modular tensor category, and if we take
the module category to be $\M \eq \C$ over $\X \eq \CbC$, i.e.\ the completion of 
\CtC\ with respect to direct sums, then $\alpha^\pm$ are monoidal equivalences.

The proof proceeds in two steps. We first compute $\intEnd (\one)$, the internal 
End of the monoidal unit $\one \iN \obj(\C)$. Assuming that \C\ is premodular 
(see Section \ref{sec_dirsums}), the internal End of any object of \M\ exists 
and is an algebra in \CbC, and the category of modules over this algebra is 
equivalent, as a module category over \CbC, to \C. 
We also show, in Section \ref{sec:frob}, that the algebra $\intEnd (\one)$
has the structure of a (symmetric special) Frobenius algebra. If \C\ is even 
a modular tensor category, i.e.\ if it is premodular and the braiding obeys a 
certain non-degeneracy condition, we show that $\intEnd(\one)$ is Azumaya.
An algebra in a braided monoidal category is an Azumaya algebra 
iff the braided induction functors $\alpha^\pm$ are monoidal equivalences,
see section \ref{sec:azumaya}. In Lemma \ref{ssFA-is-Azu} we establish an 
alternative criterion for an algebra to be Azumaya, and in Proposition \ref{A-is-Azu}
we verify that this criterion is met for $\intEnd (\one)$. We also show that the 
composition $(\alpha^-)^{-1} \cir \alpha^+$ is naturally equivalent to the 
functor that permutes the two factors in \CtC\ (Lemma \ref{lem_Gamma_M} 
and Proposition \ref{prop-Zij,kl}).

\medskip

The latter result has an application to conformal field theory. Let \V\ be a 
vertex operator algebra satisfying the conditions of Theorem 2.1 of \cite{huan21}, 
according to which the category \C\ of representations of \V\ is a modular tensor 
category. Now the product $\V \,{\otimes_\complex}\, \V$ satisfies the conditions 
of that theorem as well, and its representations form the category \CbC. Further,
by \cite{fuRs4,fjfrs}, a full conformal field theory with left and right chiral symmetry
given by the product $\V \,{\otimes_\complex}\,\V$ is described by a (symmetric 
special) Frobenius algebra in \CbC. If this Frobenius algebra is Azumaya, then
the modular invariant torus partition function of the full conformal field theory
is of permutation type. Our results imply that the algebra $\intEnd (\one)$ 
realises the conformal field theory associated to the transposition symmetry of 
the two factors in $\V \,{\otimes_\complex}\, \V$. This is non-trivial because a 
modular invariant bilinear combination of characters is not automatically the 
partition function of a physical conformal field theory with a consistent
collection of correlators. As an immediate corollary, we can also obtain algebras 
that implement the modular invariants for $\C^{\boxtimes N}$ corresponding to 
arbitrary transpositions in $\mathfrak S_N$. Since any permutation in $\mathfrak S_N$ 
can be written as a product of transpositions, and since the tensor product of 
algebras implements, by Proposition 5.3 of \cite{fuRs4}, the product of modular
invariants, our results imply that all permutation modular invariants in
rational conformal field theory are physical.

{}From this point of view another question is natural: What is the chiral 
conformal field theory whose chiral symmetry consists of the fixed algebra 
$(\V \,{\otimes_\complex}\, \V)^{\mathfrak S_2}$ under the action of the permutation 
group $\mathfrak S_2$ on $\V \,{\otimes_\complex}\, \V$\,? Theories for which the
chiral symmetry is given by a fixed algebra are called orbifold models.
On the category theoretical side, the orbifold by a symmetry group $G$ can be 
obtained by starting from a $G$-equivariant fusion category \cite{tura10}, from 
which one then calculates the orbifold category $\X/G$ \cite{kirI17,muge13}. Part 
of the data of a $G$-equivariant fusion category consists of a list of categories 
$\X_g$, indexed by elements $g\iN G$, together with bifunctors 
$\X_g \,{\times}\, \X_h \To \X_{gh}$ subject to certain compatibility 
conditions. From the analogy with conformal field theory one expects that an 
$N$-fold product $\C^{\boxtimes N}$ can be embedded into an 
$\mathfrak S_N$-equivariant category \X\ as its neutral component,
$\X_e \eq \C^{\boxtimes N}$. (This can possibly be done in different ways, 
which are indexed by elements of $H^2(\mathfrak S_N,U(1))$.) 
Our construction provides some of the necessary ingredients for obtaining an
$\mathfrak S_2$-equivariant fusion category with $\X_e \eq \CbC$ and 
$\X_\pi \eq \C$, namely, along with $\X_e \,{\times}\, \X_e \To \X_e$, also
the bifunctors $\X_\pi \,{\times}\, \X_e \To \X_\pi$ 
and $\X_e \,{\times}\, \X_\pi \To \X_\pi$
with mixed associators that obey all mixed pentagon relations.

\medskip

This paper is organised as follows. In section \ref{sec:familiy} we obtain 
the family of module category structures of \C\ over \CtC, and 
in section \ref{sec:braid} we discuss the braided induction associated to 
these module categories. For passing from \CtC\ to 
\CbC\ we need to verify that the construction is compatible 
with direct sums; this is done in section \ref{sec_dirsums}. In section 
\ref{sec_alg} we compute the internal End of $\one \iN \obj(\C)$, and in 
section \ref{sec:azumaya} we show that it is an Azumaya algebra.
In section \ref{sec:final} we conclude with some remarks and observations 
regarding our results.


\section{A family of module categories}\label{sec:familiy}

We start by recalling some basic notions that will be needed.

A (right) {\em module category\/} $(\M,\boxtimes,\psi,r)$ over a 
monoidal category $\X\eq(\X,\otimes,\one,a,\lambda,\rho)$ is \cite{ostr}\,%
 \footnote{~In \cite{ostr} it is assumed that the categories \X\ and \M\ are
abelian and that the tensor product bifunctor is exact. These assumptions 
will actually be met in the situation studied further below, in which \X\ is
premodular, but they are not needed here.} 
a category \M\ together with a bifunctor
 $\boxtimes{:}\ \M\,{\times}\,\X\,{\to}\,\M$
and natural associativity and unit isomorphisms
$\psi_{M,X,Y}{:}\ M\boti(X\oti Y)\,{\to}\,(M\boti X)\boti Y$ and 
$r_M{:}\ M\boti\one\,{\to}\,M$. They satisfy the pentagon identity
  \be
  \psi_{M\boxtimes X,Y,Z} \circ \psi_{M,X,Y\otimes Z}
  = (\psi_{M,X,Y} \boti \id_Z) \circ \psi_{M,X\otimes Y,Z} \circ (\id_M \boti a_{X,Y,Z})
  \label{modulecat-pentagon}\ee
for all $M\iN\obj(\M)$ and $X,Y,Z\iN\obj(\X)$ and the triangle identity
  \be
  \id_M \boti \lambda_X = (r_M \boti \id_X) \circ \psi_{M,\one,X}
  \ee
for all $M\iN\obj(\M)$ and $X\iN\obj(\X)$.

If $(\M,\boxtimes_\M,\psi^\M,r^\M)$ and $(\N,\boxtimes_\N,\psi^\N,r^\N)$ 
are module categories over the same monoidal category \X, then a (strict)\,%
 \footnote{~Following \cite{ostr}, by strictness we mean that the $\gamma^{}_{M,X}$
 are {\em iso\/}morphisms. In the present paper all module functors are strict, 
 so we drop this qualifier in the sequel.}
{\em module functor\/} $(G,\gamma)$ from \M\ to \N\ is \cite[Def.\,2.7]{ostr} a 
functor $G{:}~\M\To\N$ together with a collection of natural isomorphisms
  \be
  \gamma^{}_{M,X}:\quad G(M{\boxtimes_\M}X) \to G(M) \,{\boxtimes_\N}\, X
  \ee
for $M\iN\obj(\M)$ and $X\iN\obj(\X)$ that obey the following associativity
and unit constraints: first, for all $M\iN\obj(\M)$ and $X,Y\iN\obj(\X)$ one has
  \be
  (\gamma^{}_{M,X} \boxtimes_\N  \id_Y) \circ \gamma^{}_{M\boxtimes_\M X,Y}
  \circ G(\psi^\M_{M,X,Y}) = \psi^\N_{G(M),X,Y} \circ \gamma^{}_{M,X\otimes Y}
  \label{modfunct_ass}
  \ee
as morphisms $ G(M\,{\boxtimes_\M}\,(X\Oti Y)) \To (G(M) \,{\boxtimes_\N}\, X) 
\,{\boxtimes_\N}\, Y$; and second, for all $M\iN\obj(\M)$ one has
  \be
  G(r_M^\M)= r_{G(M)}^\N\cir \gamma_{M,\one}
  \ee
as morphisms $G(M\,{\boxtimes_\M}\,\one) \To G(M)$.

Given two module functors $(G,\gamma)$ and $(G',\gamma')$ from $\M$ to $\N$, a
natural transformation $\eta\colon (G,\gamma)
     $\linebreak[0]$
{\to}\,(G',\gamma')$ of module
functors is a natural transformation $\eta\colon G \To G'$ that obeys
  \be
  (\eta_M^{} \,{\boxtimes_\N}\, \id_X^{}) \circ \gamma_{M,X}^{} = \gamma'_{M,X} \circ 
  \eta_{M\boxtimes_\M X}^{}
  \ee
for all $M \iN \M$ and $X \iN \X$.
Two module categories \M\ and \N\ over \X\ are called {\em equivalent\/} iff there
exists a module functor $(F,\varphi)\colon \M \To \N$ such that $F$ is an equivalence
of categories. One can check that this is the same condition as demanding the 
existence of a module functor $(G,\gamma)\colon \N \To \M$ 
such that $(F,\varphi) \cir (G,\gamma)$ and $(G,\gamma) \cir (F,\varphi)$ are 
naturally isomorphic to $(\id_\N,\id)$ and $(\id_\M,\id)$, respectively.

\medskip

Invoking coherence, from here on we adhere to the following

\begin{conv}
Unless noted otherwise, a monoidal category is assumed to be strict, i.e.\ to
have trivial associativity and unit isomorphisms $a_{X,Y,Z}$
and $\lambda_X,\,\rho_X$. 
\\
Accordingly we will write
$\X\eq(\X,\otimes,\one)$ in place of $\X\eq(\X,\otimes,\one,a,\lambda,\rho)$.
\end{conv}

For any two categories \C\ and $\mathcal D$ we denote by $\C{\times}\,\mathcal D$ 
the category whose objects are ordered pairs $U\,{\times}\,V$ of objects of \C\ 
and $\mathcal D$ and whose morphism sets are Cartesian products 
$\Hom_{\C\times\mathcal D}(U\,{\times}\,V,U'\,{\times}\,V') \eq 
\Hom_\C(U,V) \,{\times}\, \Hom_{\mathcal D}(U',V')$ of morphism sets of \C\ and 
$\mathcal D$. If \C\ and $\mathcal D$ are monoidal, then so is $\C{\times}\,\mathcal D$, 
with tensor product $\otimes$ defined component-wise, i.e.\ as  $(U\,{\times}\,V)
     $\linebreak[0]$
{\otimes}\,
(U'\,{\times}\,V') \eq (U\,{\otimes_\C}\,U')\,{\times}\,(V\,{\otimes_{\mathcal D}}\,V')$
on objects and analogously on morphisms and on associators.

For $\C\eq\mathcal D$ this product of categories turns out to have a family of
module categories. These provide a categorification of the corresponding
structure for rings (see the introduction). Accordingly, to define them
we first introduce a bifunctor $\boxtimes{:}~ \C \times (\CtC)\To \C$ by
  \be
  M \boti (U\,{\times}\,V) := M \oti U \oti V
  \qquad \text{and} \qquad
  m \boti (f\,{\times}\,g) := m \oti f \oti g
  \label{def-boxtimes}\ee
for objects $M$ in \C\ and $U\,{\times}\,V$ in \CtC\ 
and morphisms $m{:}~ M \To M'$ and $f \,{\times}\, g{:}~ U\,{\times}\,V \To
U'\,{\times}\,V'$. 
In order to obtain the structure of a module category on \C\ we need more data,
and to this end we have to require additional structure for \C, analogous to the 
commutativity requirement in the case of rings: \C\ must be braided. We denote a 
braided monoidal category by $(\C,\otimes,\one,c)$, where $c$ is the braiding, i.e.\ a 
family of isomorphisms $c_{U,V}{:}~ U \oti V \To V \oti U$ satisfying the conditions
of \cite[Def.\,2.1]{joSt6}. We will use the shorthand
  \be
  D_{U,V} := c_{V,U}^{} \cir c_{U,V}^{}:\quad U\oti V \To U\oti V 
  \label{def-D_UV}\ee
for the two-fold braiding.

With these preparations, we can formulate

\begin{theorem}\label{thm_onlybraided}
Let $(\C,\otimes,\one,c)$ be a braided monoidal category.
For every $n\iN\zet$, $\C^{(n)} \equiv (\C,\boxtimes,\psi^{(n)},r)$ with 
$\boxtimes$ as in \erf{def-boxtimes}, $r_M \,{:=}\, \id_M$ and
  \be
  \psi^{(n)}_{M,U_1\times V_1,U_2\times V_2}
  := \Big[ (D_{M\otimes U_1\otimes V_1,U_2})^{-n}_{} \circ 
  \big(\id_{M\otimes U_1} 
  \oti c^{}_{U_2,V_1}\big) \circ 
  \big[(D_{M\otimes U_1,U_2})^{n}_{} \oti \id_{V_1}\big] \Big]
  \otimes \id_{V_2} 
  \label{psidef}\ee
is a (right) module category over \CtC.
\end{theorem}

\begin{figure}[tb]\begin{center}
  \bp(340,291)
  \put(0,135)   {$ \psi^{(0)} ~= $}
  \put(50,80) { {\bp(0,0) \scalebox{.3}{\bild{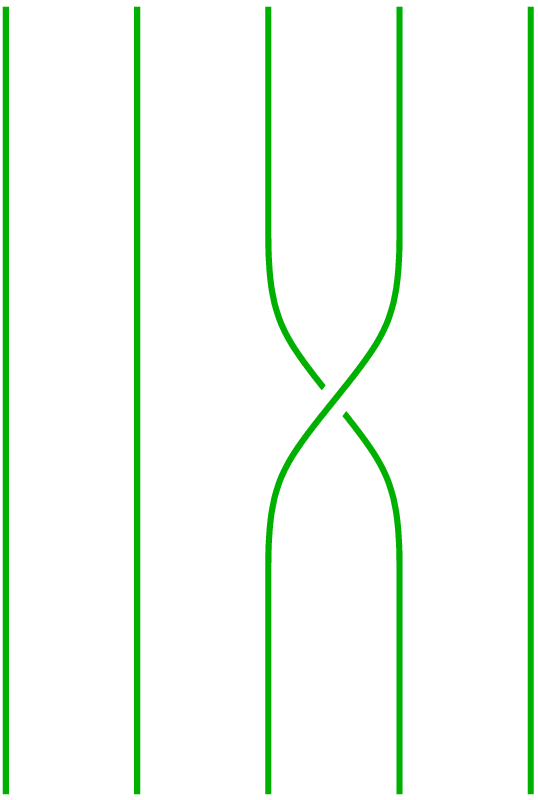}} \ep}
  \put(-4.2,-9)  {\scriptsize $ M $}
  \put(15.1,-9)  {\scriptsize $ U_1 $}
  \put(34.1,-9)  {\scriptsize $ U_2 $}
  \put(53.4,-9)  {\scriptsize $ V_1 $}
  \put(72.9,-9)  {\scriptsize $ V_2 $}
  \put(-3.2,120) {\scriptsize $ M $}
  \put(15.1,120) {\scriptsize $ U_1 $}
  \put(34.1,120) {\scriptsize $ V_1 $}
  \put(53.4,120) {\scriptsize $ U_2 $}
  \put(72.7,120) {\scriptsize $ V_2 $}
   }
  \put(200,135) {$ \psi^{(-1)} ~= $}
  \put(255,0) { {\bp(0,0) \scalebox{.3}{\bild{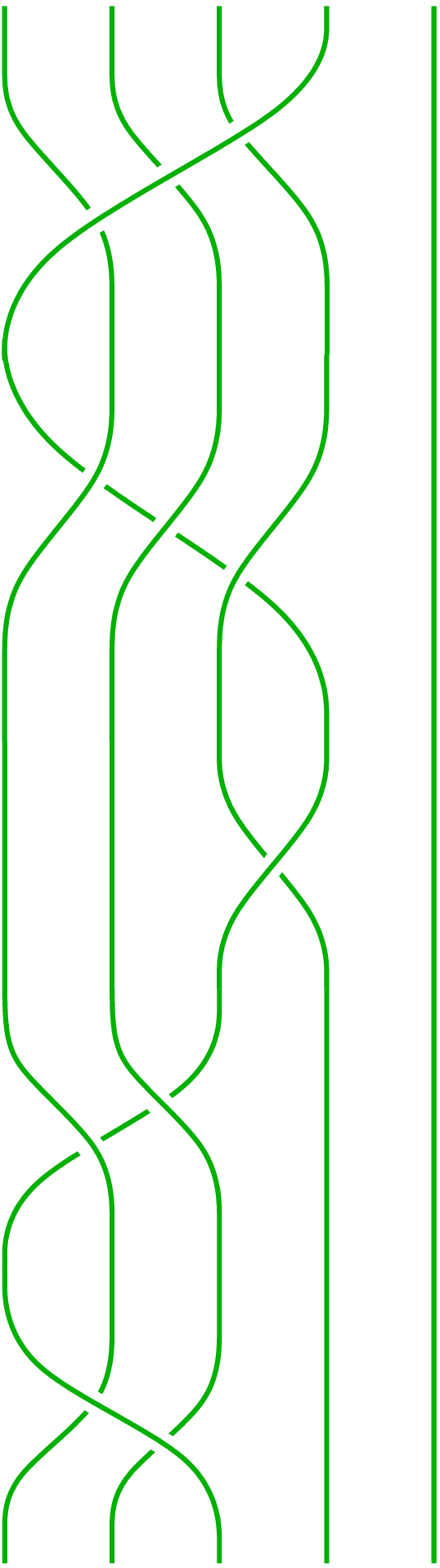}} \ep}
  \put(-4.2,-9)  {\scriptsize $ M $}
  \put(15.1,-9)  {\scriptsize $ U_1 $}
  \put(34.1,-9)  {\scriptsize $ U_2 $}
  \put(53.4,-9)  {\scriptsize $ V_1 $}
  \put(72.9,-9)  {\scriptsize $ V_2 $}
  \put(-3.2,281) {\scriptsize $ M $}
  \put(15.1,281) {\scriptsize $ U_1 $}
  \put(34.1,281) {\scriptsize $ V_1 $}
  \put(53.4,281) {\scriptsize $ U_2 $}
  \put(72.7,281) {\scriptsize $ V_2 $}
   }
  \ep
     \end{center}\caption{The morphism \erf{psidef} for $n\eq0$ and $n\eq{-}1$,
using the graphical calculus for morphisms in braided monoidal 
categories (see e.g.\ \cite{BAki,fuRs4,joSt6}).}\label{f1}
     \end{figure}

\noindent
(For the pictorial representation of $\psi^{(n)}$ in the cases $n \eq 0$ and $n \eq {-}1$ 
see figure \ref{f1})

\begin{figure}[tb]\begin{center}
  \bp(350,236)
  \put(0,0) { {\bp(0,0) \scalebox{.3}{\bild{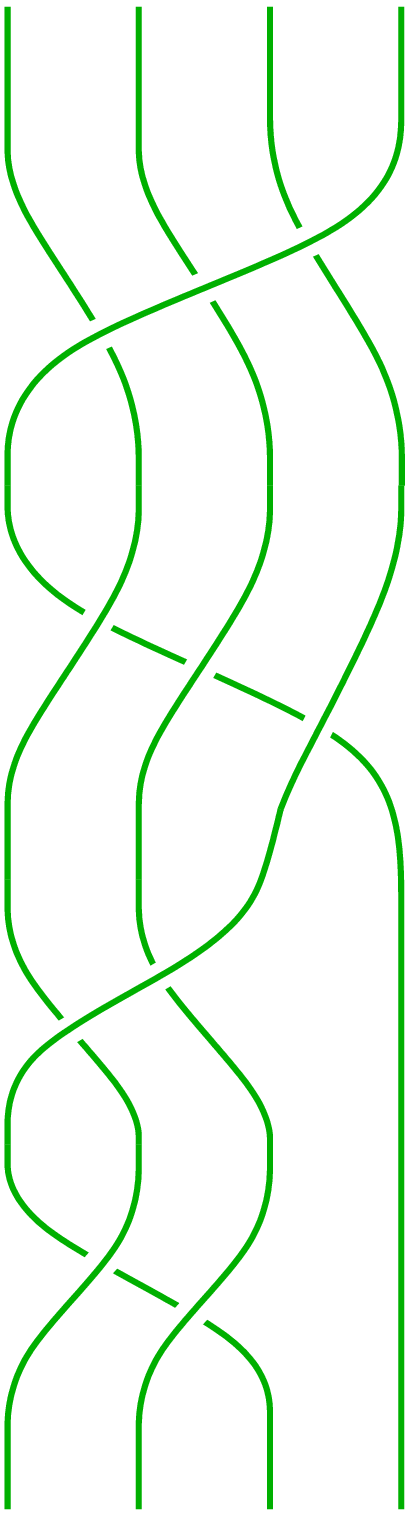}} \ep}
  \put(-4.2,-9)  {\scriptsize $ M $}
  \put(15.4,-9)  {\scriptsize $ U_1 $}
  \put(34.1,-9)  {\scriptsize $ U_2 $}
  \put(53.8,-9)  {\scriptsize $ U_3 $}
  \put(-2.6,222) {\scriptsize $ M $}
  \put(15.4,222) {\scriptsize $ U_1 $}
  \put(34.1,222) {\scriptsize $ U_2 $}
  \put(53.8,222) {\scriptsize $ U_3 $}
   }
  \put(97,106) {$ = $}
  \put(128,0) { {\bp(0,0) \scalebox{.3}{\bild{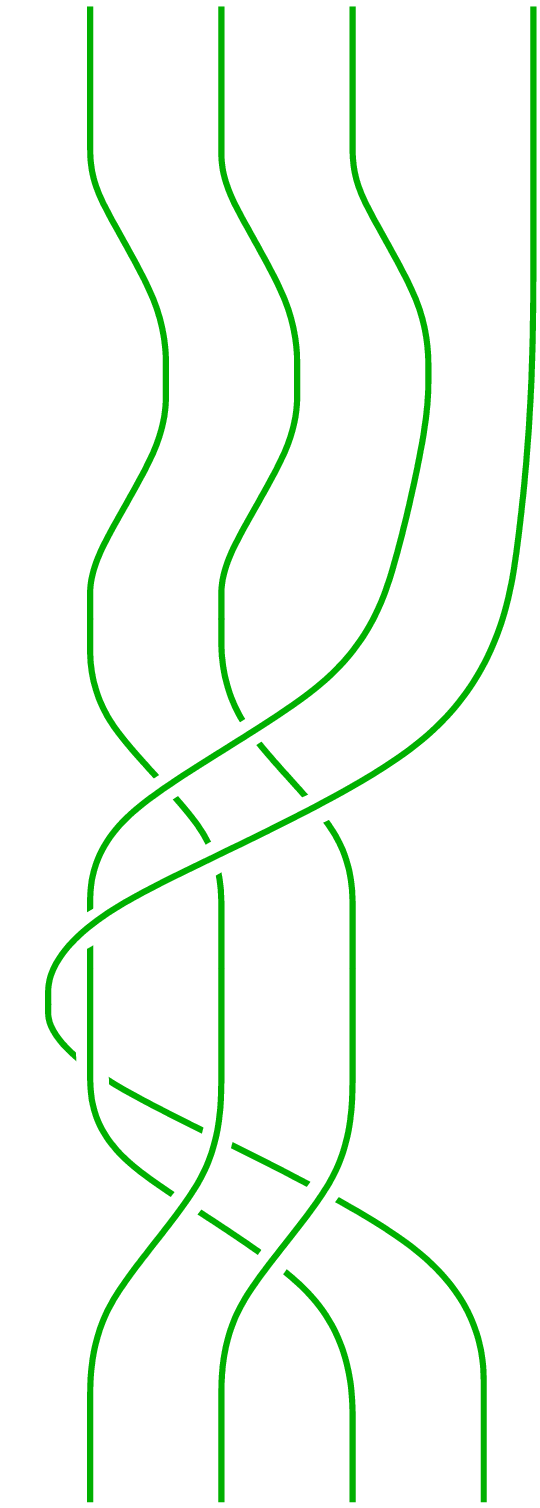}} \ep}
  \put(8.2,-9)   {\scriptsize $ M $}
  \put(28.4,-9)  {\scriptsize $ U_1 $}
  \put(47.1,-9)  {\scriptsize $ U_2 $}
  \put(66.8,-9)  {\scriptsize $ U_3 $}
  \put(8.5,222)  {\scriptsize $ M $}
  \put(28.4,222) {\scriptsize $ U_1 $}
  \put(46.9,222) {\scriptsize $ U_2 $}
  \put(72.5,222) {\scriptsize $ U_3 $}
   }
  \put(239,106) {$ = $}
  \put(275,0) { {\bp(0,0) \scalebox{.3}{\bild{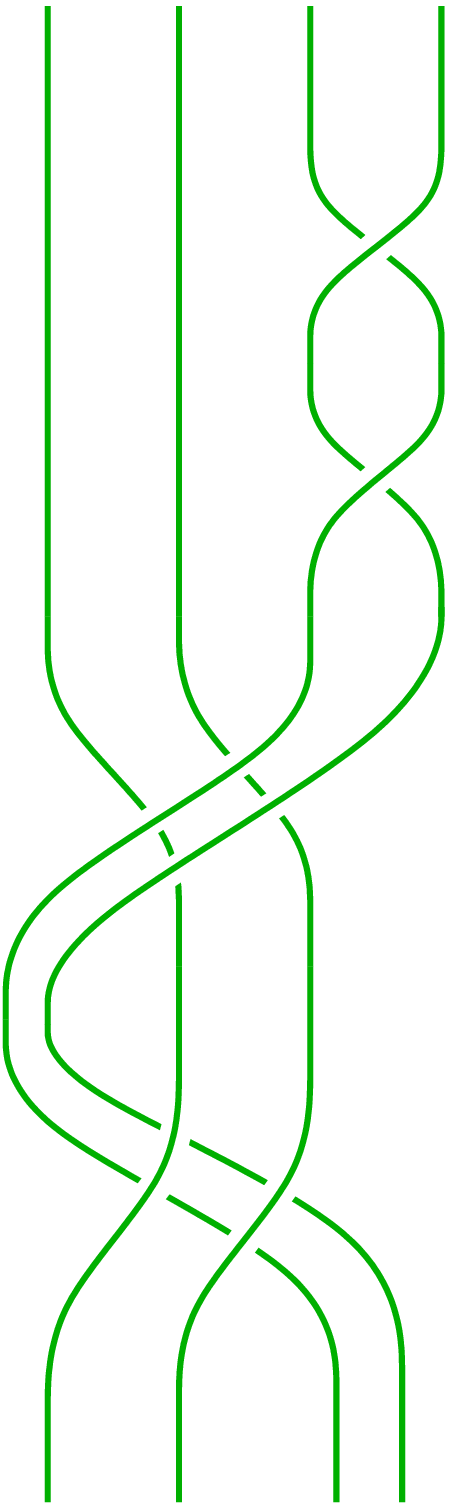}} \ep}
  \put(2.9,-9)   {\scriptsize $ M $}
  \put(23.1,-9)  {\scriptsize $ U_1 $}
  \put(44.2,-9)  {\scriptsize $ U_2 $}
  \put(55.1,-9)  {\scriptsize $ U_3 $}
  \put(3.6,222)  {\scriptsize $ M $}
  \put(22.4,222) {\scriptsize $ U_1 $}
  \put(41.1,222) {\scriptsize $ U_2 $}
  \put(59.2,222) {\scriptsize $ U_3 $}
   }
  \ep
     \end{center}\caption{Proof of the equality \erf{XXX}.}\label{f2}
     \end{figure}

\begin{proof}
Henceforth we will suppress the tensor product of \C\ on objects, e.g.\ we 
will write $UV$ instead of $U\oti V$. 
Naturality of $\psi^{(n)}$ follows from naturality of the braiding of \C, 
and the $\psi^{(n)}$ are clearly isomorphisms. The triangle identity of
$\C^{(n)}$ holds trivially. It remains to check the pentagon 
identity. Substituting the definitions, the left hand side of
\erf{modulecat-pentagon}
with $X\eq U_1\,{\times}\,V_1,~ Y\eq U_2\,{\times}\,V_2$ and 
$Z\eq U_3\,{\times}\,V_3$ reads
  \be\bearl 
  \displaystyle
  \Big\{\big[ D_{M U_1 V_1 U_2 V_2, U_3}^{-n} \circ
  \big( \id_{M U_1 V_1 U_2} \oti c^{}_{U_3,V_2} \big) \circ
  \big( D_{M U_1 V_1 U_2 , U_3}^{n} \oti \id_{V_2} \big) \big] \otimes \id_{V_3} 
  \Big\}
  \\[.5em]
  \displaystyle \qquad\quad
  \circ\, \Big\{
  \big[ D_{M U_1 V_1, U_2 U_3}^{-n} \circ
  \big( \id_{M U_1 } \oti c^{}_{U_2 U_3,V_1} \big) \circ
  \big( D_{M U_1, U_2 U_3}^{n} \oti \id_{V_1} \big) \big] \otimes \id_{V_2 V_3}
  \Big\} \,,
  \label{pentagon-lhs1}
  \eear\ee
while for the right hand side of \erf{modulecat-pentagon} one finds
  \be\bearl 
  \displaystyle
  \Big\{\big[ D_{M U_1 V_1, U_2}^{-n} \circ
  \big( \id_{M U_1} \oti c^{}_{U_2,V_1} \big) \circ
  \big( D_{M U_1, U_2}^{n} \otimes id_{V_1} \big) \big] \otimes \id_{V_2 U_3 V_3} 
  \Big\}
  \\[.5em]
  \displaystyle \qquad
  \circ\, \Big\{
  \big[ D_{M U_1 U_2 V_1 V_2 ,U_3}^{-n} \circ
  \big( \id_{M U_1 U_2} \otimes c^{}_{U_3, V_1 V_2} \big) \circ
  \big( D_{M U_1 U_2, U_3}^{n} \oti id_{V_1 V_2} \big) \big] \otimes \id_{V_3} 
  \Big\} \,.
  \label{pentagon-rhs1}
  \eear\ee
Since $D_{U,V}$ is natural in both of its arguments we can move
the morphism $D_{M U_1 U_2 V_1 V_2 ,U_3}^{-n}$ in \erf{pentagon-rhs1}
all the way to the left, where it becomes $D_{M U_1 V_1 U_2 V_2 ,U_3}^{-n}$,
so that it can then be canceled against the corresponding morphism in 
\erf{pentagon-lhs1}.
Next we rewrite the braiding $c_{U_3,V_1 V_2}$ in \erf{pentagon-rhs1} as
$c_{U_3,V_1 V_2} \eq (\id_{V_1} \oti c_{U_3, V_2}) \cir (c_{U_3,V_1} 
\oti \id_{V_2})$. Then we move the braiding $c_{U_3,V_2}$ all the way to
the left, so that we can cancel it against the corresponding braiding in 
\erf{pentagon-lhs1}. We may then also omit
the common tensor factor $\id_{V_2 V_3}$ from both sides. 
This shows that the pentagon is equivalent to
  \be\bearl 
  \displaystyle
  D_{M U_1 V_1 U_2, U_3}^n \circ D_{M U_1 V_1, U_2 U_3}^{-n}
  \circ \big( \id_{M U_1} \oti c^{}_{U_2 U_3, V_1} \big) \circ
  \big( D_{M U_1, U_2 U_3}^n \oti \id_{V_1}\big)
  \\{}\\[-.8em]
  \displaystyle
  \hspace*{3em} = \big( D_{M U_1 V_1, U_2}^{-n} \oti \id_{U_3} \big) \circ
  \big( \id_{M U_1} \oti c^{}_{U_2,V_1} \oti \id_{U_3} \big) \circ
  \big( \id_{M U_1 U_2} \oti c^{}_{U_3,V_1} \big) 
  \\[.5em]
  \displaystyle  
  \hspace*{13.9em} \circ\,
  \big(D^n_{M U_1, U_2} \oti \id_{U_3 V_1}\big) \circ
  \big(D^n_{M U_1 U_2, U_3} \oti \id_{V_1}\big) .
  \label{pentagon-both2}
  \eear\ee
Now $D_{M U_1, U_2} \oti \id_{U_3}$ and $D_{M U_1 U_2, U_3}$ commute 
owing to naturality of $D_{U,V}$. We can therefore write
$(D^n_{M U_1, U_2} \oti \id_{U_3}) \cir D^n_{M U_1 U_2, U_3}
\eq \big((D_{M U_1, U_2} \oti \id_{U_3}) \cir D_{M U_1 U_2, U_3}\big)^n$.
Further, by direct calculation using the properties of the braiding one shows that 
  \be
  (D_{M U_1, U_2} \oti \id_{U_3}) \cir D_{M U_1 U_2, U_3}
  = (\id_{M U_1} \oti D_{U_2, U_3}) \cir D_{M U_1, U_2 U_3}
  \label{XXX} \ee
(it helps to draw the corresponding pictures, see figure \ref{f2}). 
Again, the two terms on the right hand side of this equality commute, so
that we can distribute the $n$th power over both terms, after which we
can cancel $D_{M U_1, U_2 U_3}^n$ against the corresponding term on the 
left hand side of \erf{pentagon-both2}. Because of
$(c_{U_2,V_1} \oti \id_{U_3}) \cir (\id_{U_2} \cir c_{U_3,V_1})
\eq c_{U_2 U_3,V_1}$ and naturality of the braiding, the remaining braidings
now cancel from \erf{pentagon-both2}. Bringing all terms involving negative 
powers to the opposite side of the equality and noting that the various 
factors  commute thus allows us to rewrite \erf{pentagon-both2} as
  \be
  \big[ (D_{M U_1 V_1, U_2} \oti \id_{U_3}\big) \circ 
  D_{M U_1 V_1 U_2, U_3} \big]^n
  =
  \big[( \id_{M U_1 V_1} \oti D_{U_2, U_3})
  \circ D_{M U_1 V_1, U_2 U_3} \big]^n .
  \label{pentagon-both3}
  \ee
A short calculation using the properties of the braiding (again it helps to 
draw pictures) shows that the terms inside the two square brackets in 
the equality \erf{pentagon-both3} are equal, and hence that 
\erf{pentagon-both3} is indeed satisfied. 
This establishes the pentagon relation \erf{modulecat-pentagon}.
\end{proof}

We mention that for $n\eq0$ and $n\eq1$ the associativity isomorphisms involve
a single braiding:
  \be
  \psi^{(0)}_{M,U_1\times V_1,U_2\times V_2} = \id_{MU_1} \oti c^{}_{U_2,V_1}
  \oti \id_{V_2}
\ee
and
\be
  \psi^{(1)}_{M,U_1\times V_1,U_2\times V_2} = \id_{MU_1} \oti c^{-1}_{V_1,U_2}
  \oti \id_{V_2} \,.
  \label{psi0-psi1}\ee
For other values of $n$ these isomorphisms are more complicated, as 
illustrated for $n\eq{-}1$ in figure \ref{f1}.

\medskip

We also have

\begin{corollary}\label{cor_left}
For $(\C,\otimes,\one,c)$ a braided monoidal category and any
$n\iN\zet$, $(\C,\widehat\boxtimes,\widehat\psi^{(n)},l)$, with
$\widehat\boxtimes{:}~ (\CtC)\times\C \To \C$ given by
$(U\,{\times}\,V)\botih M \,{:=}\, U \oti V \oti M$ on objects
and analogously on morphisms, with $l_M \,{:=}\, \id_M$ and
with\\
 $\widehat\psi^{(n)}_{X,Y,M}\colon (X \oti Y) \boti M \To X \boti (Y \boti M)$
given by
  \be
  \widehat\psi^{(n)}_{U_1\times V_1,U_2\times V_2,M}
  := \id_{U_1} \otimes \Big[ (D_{V_1,U_2\otimes V_2\otimes M})^{n}_{}
  \circ \big( c^{-1}_{V_1,U_2} \oti \id_{V_2\otimes M} \big)
  \circ \big[ \id_{U_2} \oti (D_{V_1,V_2\otimes M})^{-n}_{} \big] \Big] ,
  \ee
is a left module category over \CtC.
\end{corollary}

\begin{proof}
The proof is completely parallel to the one of Theorem \ref{thm_onlybraided}.
When expressing morphisms via the graphical tensor calculus 
(see figure \ref{f1}), it amounts to a left-right reflection of the pictures
for the proof of Theorem \ref{thm_onlybraided}.
\end{proof}

In the sequel we concentrate on the right module category structure of \C.
Obviously, all statements have an analogue for the left module category structure.

\bigskip

Additional structure on the braided monoidal category
$\C \,{\equiv}\, (\C,\otimes,\one,c)$ can be used to relate the different
module category structures described above.
According to \cite[Def.\,6.1]{joSt6} a {\em twist\/} (or balancing) on 
\C\ is a natural family $\theta$ of isomorphisms $\theta_U{:}~ U \To U$ such
that $\theta_\one \eq \id_\one$ and
  \be
  D_{U,V} = \theta_{UV}^{} \circ (\theta_U^{-1} \oti \theta_V^{-1}) 
  \label{tensortwist}\ee
for all $U,V \iN \obj(\C)$.

\begin{theorem}\label{thm_psi-twist}
Let $\,\C$ and $\,\C^{(n)}$ be as in Theorem \ref{thm_onlybraided}.
The following are equivalent:
\\[.3em]
{\rm (i)\phantom{ii}}~For every $m,n \iN \zet$ there exists a module functor 
$(\id_\C,\gamma){:}~ \C^{(m)} \To \C^{(n)}$.
\\[.3em]
{\rm (ii)\phantom{i}}~There is a module functor $(\id_\C,\gamma){:}~ \C^{(0)} \To \C^{(1)}$.
\\[.3em]
{\rm (iii)}~\,\C\ can be endowed with a twist.
\end{theorem}

In particular, if $\C$ has a twist, then the module categories $\C^{(n)}$ 
defined in Theorem \ref{thm_onlybraided} are mutually equivalent. 
In theorem \ref{thm_psi-twist} the statement $\mathrm{(i)}\,{\Rightarrow}\,\mathrm{(ii)}$ 
is trivial. The other statements are implied by the following two lemmas.

\begin{lemma}
Let $\theta$ be a twist on $\C$. Then
$(\id_\C,\gamma){:}~\C^{(n)} \To \C^{(n+1)}$ with
  \be
  \gamma_{M,U \times V}^{} := \big[ \theta_{M \otimes U}^{-1} \circ
  (\theta_M \oti \id_U) \big] \otimes \id_V
  \label{gammadef}\ee  
for $M,U,V\iN\obj(\C)$ is a module functor.
\end{lemma}

\begin{figure}[tb]\begin{center}
  \bp(120,100)
  \put(0,38)   {$ \gamma_{M,U \times V}^{}  ~= $}
  \put(66,0) { {\bp(0,0) \scalebox{.3}{\bild{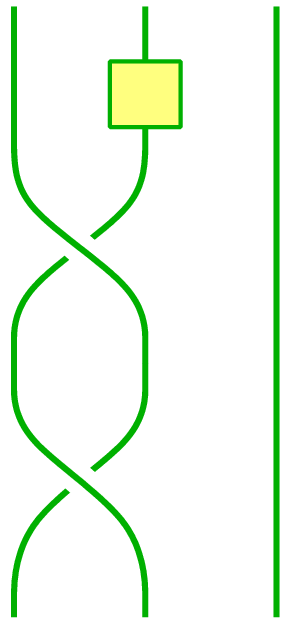}} \ep}
  \put(-2.9,-9)  {\scriptsize $ M $}
  \put(17.8,-9)  {\scriptsize $ U $}
  \put(36.8,-9)  {\scriptsize $ V $}
  \put(16.8,74.2){\scriptsize $\scriptstyle \theta^{\mbox{\tiny-}\!1} $}
  \put(-2.2,93)  {\scriptsize $ M $}
  \put(17.8,93)  {\scriptsize $ U $}
  \put(37.2,93)  {\scriptsize $ V $}
   }
  \ep
     \end{center}\caption{The morphism \erf{gammadef}.}\label{f3}
     \end{figure}

\noindent
(See figure  \ref{f3} 
 for a pictorial representation of the morphism $\gamma_{M,U \times V}^{}$.)

\begin{proof}
We must show that 
  \be
  (\gamma^{}_{M,X} \boti \id_Y) \circ \gamma^{}_{M\boxtimes X,Y} \circ 
  \psi^{(n)}_{M,X,Y} = \psi^{(n+1)}_{M,X,Y} \circ \gamma^{}_{M,X\otimes Y}
  \label{gamma-is-modulefunctor}\ee
for all $M\iN\obj(\C)$ and $X,Y\iN\obj(C{\times} C)$. 
For $X\eq U_1\,{\times}\,V_1$ and $Y\eq U_2\,{\times}\,V_2$ this can be
written as
  \be
  \begin{split}
  \hspace*{-.5em}\bearll
  \big\{ \big[ \theta^{-1}_{MU_1} \cir (\theta^{}_M \oti \id_{U_1}) \big]
  \oti \id_{V_1U_2V_2} \big\}
  \circ\big\{ \big[ \theta^{-1}_{MU_1V_1U_2} \cir (\theta^{}_{MU_1V_1} \oti
  \id_{U_2}) \big]  \oti \id_{V_2} \big\}
  \\[.5em] \quad~
  \circ \psi^{(n)}_{M,U_1\times V_1,U_2\times V_2}
  \\{}\\[-.7em] 
  = \big( D^{-1}_{MU_1V_1,U_2} \oti \id_{V_2} \big)
  \circ \psi^{(n)}_{M,U_1\times V_1,U_2\times V_2}
  \circ \big( D^{}_{MU_1,U_2} \oti \id_{V_1V_2} \big)
  \\[.5em]\quad~
\circ \big[ \theta^{-1}_{MU_1U_2} \cir (\theta^{}_M \oti \id_{U_1U_2})
  \big] \oti \id_{V_1V_2} .
  \\[-.7em]~
  \eear
  \label{gamma-is-modulefunctor-2}
\end{split}
\ee
As an illustrative example, in figure \ref{f4} we give a graphical verification in the case $n \eq {-1}$.
The general case can be proven as follows.
We express the doubled braidings in terms of twists via \erf{tensortwist}
and use the naturality of $\theta$ and of $\psi^{(n)}_{}$ to cancel various 
twist morphisms. Thereby the equality \erf{gamma-is-modulefunctor-2} can be reduced to
  \be
  \big( \theta^{-1}_{MU_1} \oti \id_{V_1U_2V_2} \big)
  \circ \psi^{(n)}_{M,U_1\times V_1,U_2\times V_2}
  = \psi^{(n)}_{M,U_1\times V_1,U_2\times V_2}
  \circ \big( \theta^{-1}_{MU_1} \oti \id_{U_2V_1V_2} \big) \,.
  \ee
That the latter equality is indeed satisfied is easily verified by
inspection of the explicit form of $\psi^{(n)}_{}$.
     \begin{figure}[!tb]\begin{center}
  \bp(430,407)
  \put(0,0) { {\bp(0,0) \scalebox{.3}{\bild{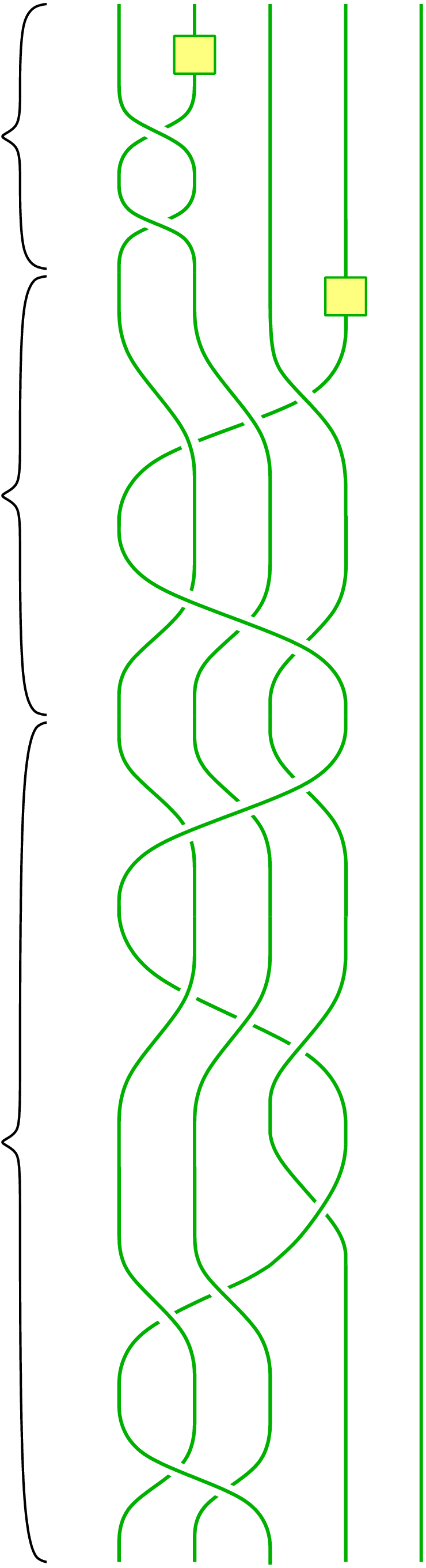}} \ep}
  \put(-12,335) {\begin{turn}{90}$ \gamma^{}_{M,X} \boti \id_Y $\end{turn}}
  \put(-11,253) {\begin{turn}{90}$ \gamma^{}_{M\boxtimes X,Y} $\end{turn}}
  \put(-17,92)  {\begin{turn}{90}$ \psi^{(-1)}_{M,X,Y} $\end{turn}}
   \put(30,0){
  \put(-5.2,-9)  {\scriptsize $ M $}
  \put(15.1,-9)  {\scriptsize $ U_1 $}
  \put(34.1,-9)  {\scriptsize $ U_2 $}
  \put(53.4,-9)  {\scriptsize $ V_1 $}
  \put(72.7,-9)  {\scriptsize $ V_2 $}
  \put(14.8,377.2){\scriptsize $\scriptstyle \theta^{\mbox{\tiny-}\!1} $}
  \put(52.8,316.5){\scriptsize $\scriptstyle \theta^{\mbox{\tiny-}\!1} $}
  \put(-4.2,397) {\scriptsize $ M $}
  \put(15.1,397) {\scriptsize $ U_1 $}
  \put(34.4,397) {\scriptsize $ V_1 $}
  \put(53.4,397) {\scriptsize $ U_2 $}
  \put(72.7,397) {\scriptsize $ V_2 $}
   }}
  \put(138,191){$ = $}
  \put(169,0) { {\bp(0,0) \scalebox{.3}{\bild{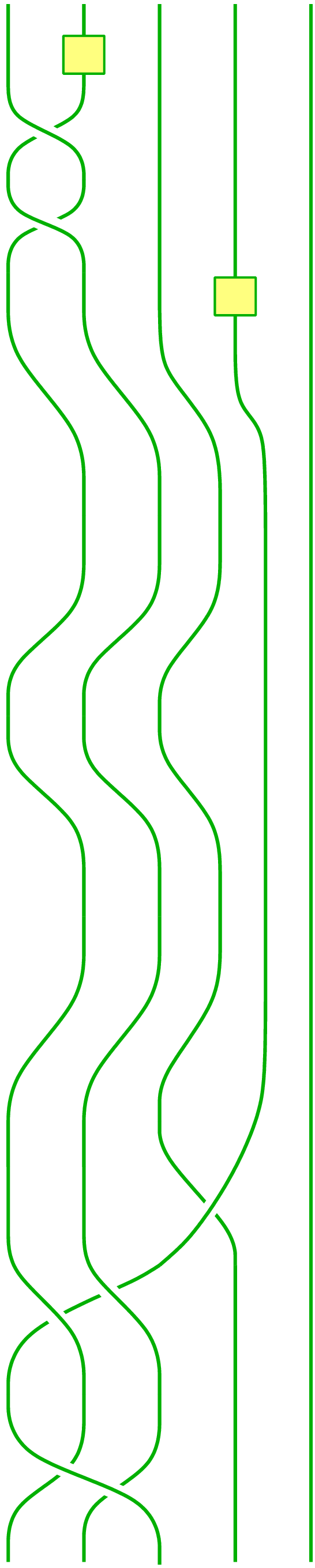}} \ep}
   \put(2.2,0){
  \put(-5.2,-9)  {\scriptsize $ M $}
  \put(15.1,-9)  {\scriptsize $ U_1 $}
  \put(34.1,-9)  {\scriptsize $ U_2 $}
  \put(53.1,-9)  {\scriptsize $ V_1 $}
  \put(72.7,-9)  {\scriptsize $ V_2 $}
  \put(14.8,377.2){\scriptsize $\scriptstyle \theta^{\mbox{\tiny-}\!1} $}
  \put(52.8,316.5){\scriptsize $\scriptstyle \theta^{\mbox{\tiny-}\!1} $}
  \put(-4.2,397) {\scriptsize $ M $}
  \put(15.1,397) {\scriptsize $ U_1 $}
  \put(34.4,397) {\scriptsize $ V_1 $}
  \put(53.4,397) {\scriptsize $ U_2 $}
  \put(72.7,397) {\scriptsize $ V_2 $}
   }}
  \put(278,191){$ = $}
  \put(311,0) { {\bp(0,0) \scalebox{.3}{\bild{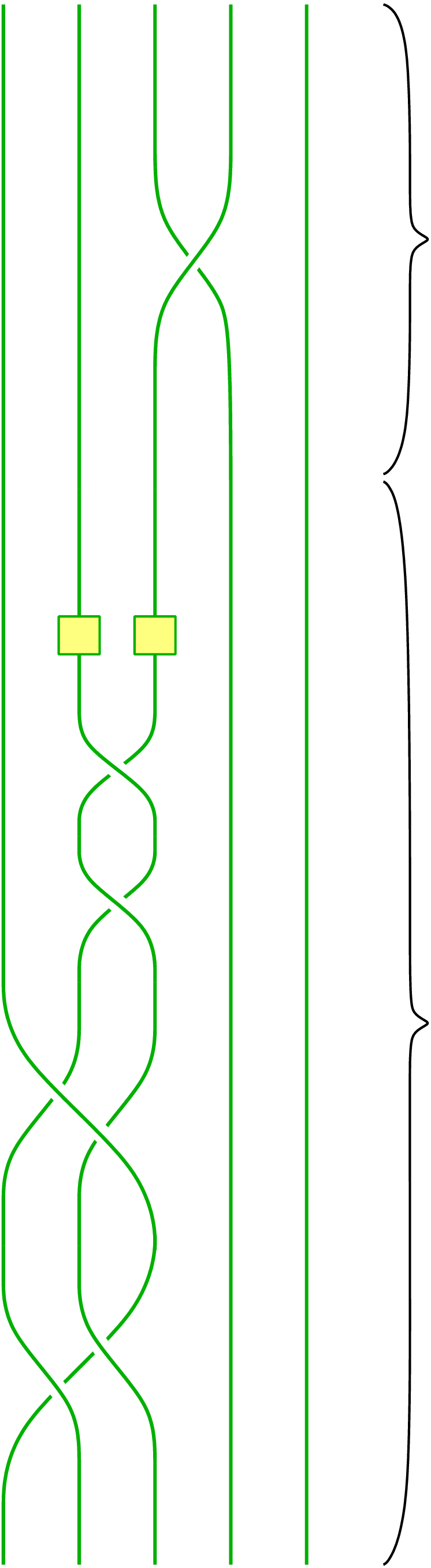}} \ep}
  \put(-5.2,-9)  {\scriptsize $ M $}
  \put(15.1,-9)  {\scriptsize $ U_1 $}
  \put(34.1,-9)  {\scriptsize $ U_2 $}
  \put(53.4,-9)  {\scriptsize $ V_1 $}
  \put(72.1,-9)  {\scriptsize $ V_2 $}
  \put(15.8,231.2){\scriptsize $\scriptstyle \theta^{\mbox{\tiny-}\!1} $}
  \put(34.8,231.2){\scriptsize $\scriptstyle \theta^{\mbox{\tiny-}\!1} $}
  \put(-3.2,397) {\scriptsize $ M $}
  \put(15.1,397) {\scriptsize $ U_1 $}
  \put(34.4,397) {\scriptsize $ V_1 $}
  \put(53.4,397) {\scriptsize $ U_2 $}
  \put(72.7,397) {\scriptsize $ V_2 $}
   }
  \put(424,350) {\begin{turn}{270}$ \psi^{(0)}_{M,X,Y} $\end{turn}}
  \put(424,155) {\begin{turn}{270}$ \gamma^{}_{M,X\otimes Y} $\end{turn}}
  \ep
     \end{center}\caption{The equality \erf{gamma-is-modulefunctor} for $n\eq{-}1$.}\label{f4}
     \end{figure}
\end{proof}

\begin{lemma}\label{lem_C0-C1}
Let $(\id_\C,\gamma){:}~\C^{(0)} \To \C^{(1)}$ be a module functor. Then
  \be
  \theta_U := \gamma_{\one,\one \times U}^{} \circ 
  (\gamma_{\one,U \times \one}^{})^{-1}_{}
  \label{def-theta}\ee
is a twist in $\C$.
\end{lemma}

\begin{proof}
By definition, $\theta_\one\eq\id_\one$, and naturality of $\theta$ follows
from naturality of $\gamma$. Thus to establish that $\theta$ is a twist, it
remains to verify the equality \erf{tensortwist}.
\\
Using the explicit form \erf{psi0-psi1} of $\psi^{(0)}_{}$ and $\psi^{(1)}_{}$,
according to \erf{gamma-is-modulefunctor} we have
  \be
  \begin{array}r
  (\gamma^{}_{M,U_1\times V_1} \oti \id_{U_2V_2})
  \circ \gamma^{}_{MU_1V_1,U_2\times V_2}
  \circ (\id_{MU_1} \oti c^{}_{U_2,V_1} \oti \id_{V_2})
  \hspace*{5em}\\{}\\[-.7em]
  = (\id_{MU_1} \oti c^{-1}_{V_1,U_2} \oti \id_{V_2})
  \circ \gamma^{}_{M,U_1U_2\times V_1V_2} \,.
  \eear\ee
When considered for $M\eq U_2 \eq V_1 \eq \one$ and $U_1\eq U$, $V_2 \eq V$,
this reduces to
  \be
  (\gamma^{}_{\one,U\times\one} \oti \id_V) \circ \gamma^{}_{U,\one\times V}
  = \gamma^{}_{\one,U\times V} \,,
  \label{220}\ee
and similarly one deduces that
  \be
  (\gamma^{}_{M,\one\times U} \oti \id_V) \circ \gamma^{}_{MU,\one\times V}
  = \gamma^{}_{M,\one\times UV} \,, \qquad
  (\gamma^{}_{\one,U\times\one} \oti \id_V) \circ \gamma^{}_{U,V\times\one}
  = \gamma^{}_{\one,UV\times\one} 
  \label{221}\ee
and
  \be
  (\gamma^{}_{\one,\one\times V} \oti \id_U) \circ \gamma^{}_{V,U\times\one}
  \circ c^{}_{U,V} = c^{-1}_{V,U} \circ \gamma^{}_{\one,U\times V} \,.
  \label{222}\ee
When combined with \erf{220} and \erf{221}, the relation \erf{222}
can be rewritten as
  \be
  \begin{array}r
  (\gamma^{}_{\one,\one\times V} \oti \id_U)
  \circ (\gamma^{-1}_{\one,V\times\one} \oti \id_U)
  \circ \gamma^{}_{\one,VU\times\one} \circ c^{}_{U,V} 
  \hspace*{12em}\\{}\\[-.7em]
  = c^{-1}_{V,U} \circ (\gamma^{}_{\one,U\times\one} \oti \id_V)
  \circ (\gamma^{-1}_{\one,\one\times U} \oti \id_V)
  \circ \gamma^{}_{\one,\one\times UV} \,.
  \eear
  \label{13}\ee
Next we note that naturality of $\gamma$ implies that
  \be
  (\id_M\oti f\oti g) \circ \gamma^{}_{M,R_1\times S_1}
  = \gamma^{}_{M,R_2\times S_2} \circ (\id_M\oti f\oti g)
  \ee
for all morphisms $f\,{\times}\,g{:}~ R_1\,{\times}\,S_1 \To
R_2\,{\times}\,S_2$ in \CtC. For $R_1\eq UV$, $R_2\eq VU$,
$M\eq S_1\eq S_2\eq\one$ and $f\eq c^{}_{U,V}$, $g\eq\id_1$, this reads
  \be
  \gamma^{}_{\one,VU\times\one} \circ c^{}_{U,V}
  = c^{}_{U,V} \circ \gamma^{}_{\one,UV\times\one} \,.
  \ee
Combining this equality with \erf{13} and with the naturality of $c$, we finally obtain
  \be
  (\theta_V^{} \oti \theta_U^{}) \circ c^{}_{U,V}
  = c^{-1}_{V,U} \circ  \theta_{UV}^{}  
  \ee
with $\theta$ as in \erf{def-theta}. As a consequence, $\theta$ satisfies
\erf{tensortwist}, as required.
\end{proof}

\begin{rem}\label{rem_G_gamma-vs-id_gamma}
In Theorem \ref{thm_psi-twist} we only consider module functors of the form 
$(\id_\C,\gamma)\colon \C^{(m)} \To \C^{(n)}$. But in fact condition (i) of Theorem 
\ref{thm_psi-twist} can be replaced by the weaker requirement that for all $m,n\iN\zet$ 
the module categories $\C^{(m)}$ and $\C^{(n)}$ are equivalent as module categories. To 
prove this, by Lemma \ref{lem_C0-C1} it is enough to show that the existence of an 
equivalence $(G,\gamma) \colon \C^{(0)} \To \C^{(1)}$ implies the existence of an 
equivalence $(\id,\tilde\gamma) \colon \C^{(0)} \To \C^{(1)}$.
\\
For $Y \iN \obj(\C)$ denote by $T_Y$ the functor $Y{\otimes}- \colon \C \To \C$. One 
checks that $(T_Y,\id) \colon \C^{(0)} \To \C^{(0)}$  is a module functor. Next one 
notes that $\gamma_{\one,U\times\one} \colon G(U) \To G(\one) \oti U$ is a natural 
isomorphism of functors from $G$ to $T_{G(\one)}$. Since $G$ is an equivalence of module 
categories, there exists a functor $(F,\varphi)$ such that $(G,\gamma) \cir (F,\varphi)$ is 
naturally isomorphic to $(\id_\C,\id)$ as a module functor. The composition 
$(G,\gamma) \cir (T_{F(\one)},\id)$ is then naturally isomorphic to 
$(\id_\C,\tilde\gamma)$ for some $\tilde\gamma$.
\end{rem}


\section{Braided induction}\label{sec:braid}

For the construction of the module category structures of $\M\eq\C$ over $\X\eq\CtC$ 
  in the previous section one must assume that \C\ is braided. From this assumption
  it follows that \X\ can be endowed with a braiding, too:\,%
 \footnote{~In fact there are several braidings, depending on whether $c$ or $c^{-1}$
  is taken in each factor. For us the choice $c \,{\times}\, c$ is the relevant one.}
    \be
    c_{U\,{\times}\,V,U'\,{\times}\,V'}^{\C\times\C} = c_{U,U'}^{} \,{\times}\, c_{V,V'}^{} \,,
    \label{c_for_CxC}\ee   
  where $c$ is the braiding of \X. Now it is known that to any module category
  over a braided category \X\ there are associated two functors, called
  braided induction, from \X\ to the category of module endofunctors.
These turn out to be of interest also in the present context.

\medskip

If $(\M,\boxtimes,\psi)$ is a module category over a braided monoidal category \X, 
then for any $X\iN\obj(\X)$ one defines \cite{ostr} two module endofunctors 
$\alpha_X^\pm$ of \M\ as follows. As functors they coincide, both acting as
  \be
  \alpha_X^\pm(M) := M \boti X \qquad{\rm and}\qquad
  \alpha_X^\pm(f) := f \boti \id_X 
  \label{def_alpha}\ee
on objects and morphisms of \M, while their module functor structures  
$\gamma^{X,\pm}$ are given by 
  \be
  \bearl
  \gamma^{X,+}_{M,Y} := \psi^{}_{M,X,Y} \circ (\id_M \boti c_{Y,X}^{}) \circ \psi^{-1}_{M,Y,X} 
  \qquad{\rm and} \\{}\\[-.8em]
  \gamma^{X,-}_{M,Y} := \psi^{}_{M,X,Y} \circ (\id_M \boti c_{X,Y}^{-1}) \circ \psi^{-1}_{M,Y,X} \,,
  \eear
  \label{gamma_X_pm}\ee
respectively.

For the module categories constructed in Theorem \ref{thm_onlybraided}, it is 
straightforward to specify natural transformations between the endofunctors $\alpha^\pm$ 
-- these will be instrumental in the proof of Proposition \ref{prop-Zij,kl} below.
Only the case $n\eq 0$ of Theorem \ref{thm_onlybraided} will be needed for this purpose.

\begin{lemma}\label{lem_Gamma_M}
Let $\C^{(0)}$ be the module category over $\CtC$ given in Theorem \ref{thm_onlybraided}.
Then for any pair of objects $U,V\iN\obj(\C)$, the collection
  \be
  \Gamma_{\!M} := [ (c_{V,M}^{} \cir c_{M,V}^{}) \oti \id_U ] \circ
  (\id_M \oti c_{U,V}^{})
  \label{Gamma_M}\ee
of morphisms in $\Hom_\C(M\oti U\oti V,M\oti V\oti U)$ furnishes a natural 
isomorphism from $\alpha^+_{U\times V}$ to $\alpha^-_{V\times U}$ as module functors.
\end{lemma}

\begin{proof}
First note that $\alpha^\pm_{U\times V}(M) \eq M\oti U\oti V$, so that the
morphisms \erf{Gamma_M} belong to the correct morphism spaces. In view of the 
definition \erf{def_alpha} of $\alpha^\pm_X$ and formulas \erf{psi0-psi1} and 
\erf{c_for_CxC} for the associator and the braiding, it is then obvious that these 
isomorphisms constitute a natural transformation from $\alpha^+_{U\times V}$ to 
$\alpha^-_{V\times U}$ as functors.
\\
It thus remains to be checked that this natural transformation is compatible
with the module functor property. This means that the equality
  \be
  \gamma^{V\times U,-}_{M,U'\times V'} \circ \Gamma_{\!M\boxtimes(U'\times V')}
  = (\Gamma_{\!M} \oti \id_{U'\otimes V'}) \circ \gamma^{U\times V,+}_{M,U'\times V'} \,,
  \label{gamma_Gamma}\ee
with morphisms $\gamma^{X,\pm}$ as in \erf{gamma_X_pm},
holds for all $M\iN\obj(\C)$ and $U'{\times}V'\iN\obj(\CtC)$. Now according 
to \erf{psi0-psi1} and \erf{c_for_CxC},
not only the $\Gamma_{\!M}$, but also the relevant morphisms 
$\gamma^{X,\pm}$ are entirely expressible through the braiding of \C:
  \be
  \bearl
  \gamma^{U\times V,+}_{M,U'\times V'} = \id_M \otimes [\,
  (\id_{U} \oti c_{U',V}^{} \oti \id_{V'}) \circ (c_{U',U}^{} \oti c_{V',V}^{})
  \circ (\id_{U'} \oti c_{U,V'}^{-1} \cir \id_{V}) \,] \,,
  \\{}\\[-.8em]
  \gamma^{V\times U,-}_{M,U'\times V'} = \id_M \otimes [\,
  (\id_{V} \oti c_{U',U}^{} \oti \id_{V'}) \circ (c_{V,U'}^{-1} \oti c_{U,V'}^{-1})
  \circ (\id_{U'} \oti c_{V,V'}^{-1} \cir \id_{U}) \,] \,.
  \eear
  \ee
It is then an easy application of the properties of the braiding
(amounting here to using relations in the braid group on five strands) to check that
the equalities \erf{gamma_Gamma} are indeed satisfied.
\end{proof}

The composition of functors endows the category $\Endfun\X\M$ of
module endofunctors with the structure of a strict monoidal category.
(Also, \X\ is a left module category over $\Endfun\X\M$.)
The functors $\alpha^\pm$ of {\em braided induction\/} (or $\alpha$-{\em induction\/}
\cite{lore,xu3,boev}) are the two functors from \X\ to $\Endfun\X\M$ that act as
  \be
  \alpha^\pm_{}:\quad \obj(\X) \ni X \,\longmapsto\,  \alpha^\pm_X
  \label{alphapm}\ee
on objects. One checks that $\alpha^\pm$ together with the associativity and unit 
constraints of \M\ constitute {\em monoidal\/} functors from $(\X,\otimes^{\rm opp})$ 
to $\Endfun\X\M$. 
The braided induction functors are of particular interest in the case that \X\ is a 
modular tensor category and a few other conditions (see \cite[Thm\,O]{fuRs12}) on the 
module category are satisfied. In that situation every module endo\-functor of \M\ can be 
obtained as a retract, in the monoidal category $\Endfun\X\M$, of one of the 
objects $\alpha_X^+\cir \alpha_Y^-\iN\Endfun\X\M$, see Theorem 4.10 of \cite{ffrs}. 
Crucial information about the module category \M\ over \X\ is then encoded in the spaces
$\Hom_{\Endfun\X\M}(\alpha_X^+,\alpha_Y^-)$ of natural transformations between 
braided induction functors.


\section{Completion with respect to direct sums}\label{sec_dirsums}

In the sequel we will be interested in the case that \C\ is {\em premodular\/}
with ground field $\complex$, that is \cite{brug2}, an abelian 
$\complex$-linear semisimple category with finitely many isomorphism classes 
of simple objects, endowed with a ribbon structure
(i.e.\ being braided monoidal and having a compatible twist and dualities, 
see e.g.\ \cite{KAss}), and with simple monoidal unit. Then 
the tensor product of \C\ is exact and thus preserves direct sums. As a consequence, 
also the bifunctor $\boxtimes$ defined as in \erf{def-boxtimes} preserves direct sums.

Braided induction is compatible with direct sums as well. Indeed,
it is straightforward to check that, for \M\ a module category over \X\ and 
$\{X_\ell\}$ any family of objects in \X\ whose direct sum $\bigoplus_\ell X_\ell$
exists, for every $M\iN\obj(\M)$ the object $\alpha^\pm_{\bigoplus_\ell\! 
X_\ell}(M)$ of \M\ satisfies the universal property of the direct sum 
$\bigoplus_\ell(\alpha^\pm_{X_\ell}(M))$.

In the case of premodular
categories the following notion of product of categories is relevant.

\begin{definition}
Let \C\ be an abelian $\complex$-linear category. Then the category $\CbC$ is defined 
as follows. $\obj(\CC)$ is obtained by completing the class of pairs 
$U\,{\times}\,V$ of objects of \C\ with respect to direct sums. The morphism sets 
of \CbC\  are tensor products $\Hom_{\CC}(U\,{\times}\,V,U'\,{\times}\,V') \eq $
$\Hom_\C(U,V)\,{\otimes_\complex}\,\Hom_\C(U',V')$ of morphism spaces of \C, extended to
direct sums.
\end{definition}

If \C\ is premodular, then so is \CbC, with tensor product (denoted by $\otimes$ again) 
and braiding (denoted by $c$ again) defined in the same way
as in the case of \CtC, extended to direct sums.
As a consequence, from Theorem \ref{thm_onlybraided} we can deduce that
for \C\ a $\complex$-linear braided monoidal category,
for any $n\iN\zet$, $\C^{(n)}\,{\equiv}\,(\C,\boxtimes,\psi^{(n)})$ with
  \be
  \boxtimes{:}\quad M \times (\bigoplus_\ell U_\ell\,{\times}\,V_\ell)
  \,\mapsto\, \bigoplus_\ell M \oti U_\ell \oti V_\ell
  \ee
and
  \be
  \bearl
  \psi^{(n)}_{M,\bigoplus_l(U_l\times V_l),\bigoplus_m(U_m\times V_m)}
  \\{}\\[-.7em] \hspace*{2.5em}\dsty
  :=\bigoplus_{l,m} \Big[ (D_{MU_lV_l,U_m})^{-n} \circ \big(\id_{MU_l} 
  \oti c^{}_{U_m,V_l}\big) \circ 
  \big[ (D_{MU_l,U_m})^{n} \oti \id_{V_l}\big] \Big] \otimes \id_{V_m} 
  \eear
  \ee
is a (right) module category over \CbC.

Similarly, from  Theorem \ref{thm_psi-twist} and Remark \ref{rem_G_gamma-vs-id_gamma}
we learn that the following statements about $(\C,\boxtimes,\psi^{(n)})$ are equivalent:
\\[.1em]
{\rm (i)\phantom{ii}}~For every $m,n \iN \zet$ there exists a module functor 
$(\id_\C,\gamma){:}~ \C^{(m)} \To \C^{(n)}$.
\\[.1em]
{\rm (ii)\phantom{i}}~\,\C\ can be endowed with a twist.
\\[.17em]
{\rm (i$'$)}~\,\,For every $m,n \iN \zet$ the module categories $\C^{(m)}$ and $\C^{(n)}$ 
are equivalent.

\begin{conv} \label{C-premodular}
In the sequel, unless stated otherwise, the category \C\ will be assumed to be
a premodular category. We choose a 
set $\{ U_i \,|\, i \iN I \}$ of representatives of the simple objects of \C,
and we set $\NN ijk\,{:=}\,\dimc(\Hom_\C(U_i\oti U_j,U_k))\,$.
\end{conv}


\section{A family of algebras in { \CbC}}\label{sec_alg}

Recall from \cite[Thm.\,1]{ostr} that any semisimple indecomposable (right) module category 
\M\ over a semisimple rigid monoidal category \X\ with finitely many isomorphism classes 
of simple objects and with simple unit object is equivalent to the category of
(left) $A$-modules for some algebra $A$ in \X. Such an algebra can be obtained as the 
internal End $\intEnd (M)$ of a non-zero object
$M$ in \M; different objects of $\M$ lead to Morita equivalent algebras.

First recall from \cite[Sec.\,3.2]{ostr} that for two objects $M_1$ and $M_2$ of \M their internal Hom 
$\underline\Hom (M_1,M_2)$ is
an object of \X\ that represents the functor
\be
\begin{array}{ll}
\X&\to\vectk\\
X&\mapsto \Hom_\M(M_1\boti X,M_2).
\end {array}
\ee
The internal Hom always exists, when \X\ and \M\ are semisimple and \X\ has only finitely many isomorphism classes of simple objects.
For $M=M_1=M_2$ we write $\intEnd(M):=\underline\Hom(M,M)$ and call this object the internal End of $M$.
There is a canonical evaluation morphism $e_{M_1,M_2}:M_1\boti\underline\Hom(M_1,M_2)\to M_2$ in \M that induces
an associative composition morphism 
\be
\underline\Hom(M_1,M_2)\oti\underline\Hom(M_2,M_3)\to\underline\Hom(M_1,M_3).
\ee
This morphism gives the structure of an associative algebra on the internal End of an object $M$ in \M.

In the situation at hand there is a distinguished object in \M, namely the tensor unit 
$\unit$. We are interested in the algebra structure of the internal End of $\unit$. 
The multiplication morphism depends on which of the module category
structures $\C^{(n)}$ of section \ref{sec_dirsums} we choose; accordingly we 
denote the internal End of $\unit$ by $\intEnd^{(n)}(\unit)$.
In the sequel we will again use the graphical calculus for morphisms, now for the case of
premodular categories.
We depict the morphisms in dual bases $\{f_\alpha\}$ of
$\Hom_\C(U_i\oti U_j,U_k)$ and $\{\bar f_\alpha\}$ of $\Hom_\C(U_k,U_i\oti U_j)$ as 
  \be\label{f-fbar-basis}
\raisebox{-30pt}{
  \bp(170,71)
  \put(-40,31)   {$ f_\alpha ~= $}
  \put(0,0)   {\bp(0,0) \scalebox{.4}{\bild{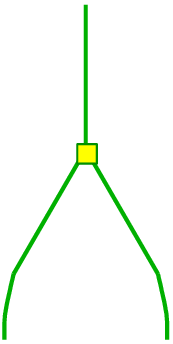}} \ep}
  \put(-2.5,-8)  {\scriptsize $ U_i $}
  \put(12.5,68.2){\scriptsize $ U_k $}
  \put(28.2,-8)  {\scriptsize $ U_j $}
  \put(20.1,35)  {\scriptsize $ f_\alpha $}
  \put(70,31)    {and $\qquad \bar f_{\alpha} ~= $}
  \put(165,0) {\bp(0,0) \scalebox{.4}{\bild{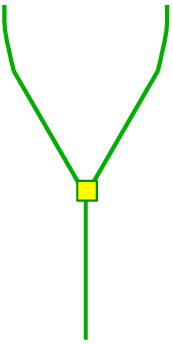}} \ep
    \put(-2.5,68.2) {\scriptsize $ U_i $}
    \put(12.5,-8)   {\scriptsize $ U_k $}
    \put(28.6,68.2) {\scriptsize $ U_j $}
    \put(20.1,27)   {\scriptsize $ \bar f_\alpha $} }
  \ep}
  \ee
~\\

\begin{theorem}\label{thm:End1=A}
Let \C and $\C^{(n)}$ be as in section \ref{sec_dirsums}. As an object of \CbC,
  \be
  \intEnd^{(n)}(\unit) \,\cong\, \bigoplus_{i\iN\calI}\,U_i\Vee\times U_i =: A \,.
  \ee
The multiplication of the algebra $\intEnd^{(n)}(\unit)=(A,m^{(n)},\eta)$ is given by
  \be\label{multiplication-morphism}
\raisebox{-54pt}{
  \bp(370,135)
\setlength{\unitlength}{1.2pt}
  \put(24,49)     {$ m^{(n)} ~=~ \dsty\bigoplus_{i,j,k\iN\I}\, \sum_{\alpha=1}^{\NN ijk}$}
  \put(127,0) {\bp(0,0) \scalebox{.48}{\bild{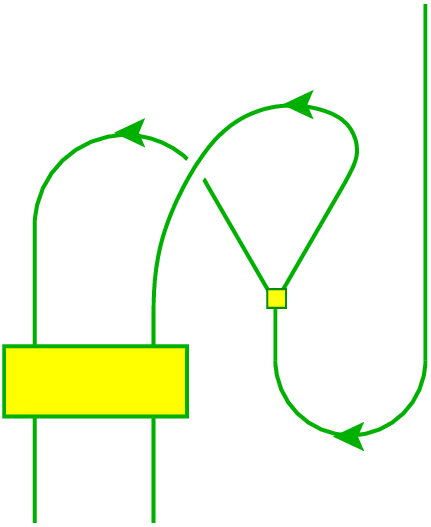}} \ep
    \put(.5,-8)   {\scriptsize $ U_i^\vee $}
    \put(23.2,-8)   {\scriptsize $ U_j^\vee $}
    \put(76.5,109)  {\scriptsize $ U_k^\vee $}
    \put(6.2,28.3)  {\scriptsize $ D_{U_i^\vee\!,U_j^\vee}^n$}
    \put(56.1,42.1) {\scriptsize $ \bar f_\alpha $} }
  \put(223,49)  {$ \mortimes $}
    \put(245,0) {\bp(0,0) \scalebox{.48}{\bild{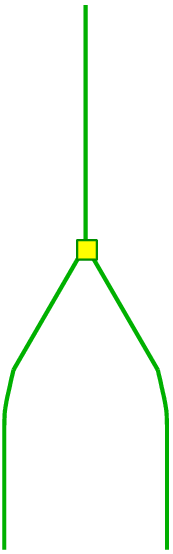}} \ep
    \put(-2.8,-8)   {\scriptsize $ U_i $}
    \put(12.5,110)  {\scriptsize $ U_k $}
    \put(28.2,-8)   {\scriptsize $ U_j $}
    \put(20.1,57)   {\scriptsize $ f_\alpha $} }
  \ep}
  \ee
~\\
The unit morphism is $\eta \eq e_{\one\times\one\prec A}$, the embedding 
morphism of $\,\one\,{\times}\,\one$ as a subobject of $A$.
\end{theorem}

\noindent
For $n\eq 0$ this algebra has been studied, in different contexts, in 
\cite{muge8,ffrs,koRu}. Observe that the product $m^{(n)}$ is independent of 
the choice of bases of the morphism spaces $\Hom_\C(U_i\oti U_j, U_k)$. 

\begin{proof}
We perform the following steps:\\
(1)~Determine the structure of an internal End on the object $A$.\\
(2)~Find the evaluation morphism \cite{ostr} $\ev{:}~\unit\boti A \To \unit$.\\
(3)~Compute the multiplication morphism on $A$.
\\
Only the last step involves the different structures of module categories on \C.

\medskip\noindent
(1)~ We consider the functors
  \be
  \begin{array}{lll} 
  & T{:}\quad \CC \To \calC\,,\quad&
  \bigoplus_{\ell}V_\ell\times W_\ell\mapsto \bigoplus_{\ell}V_\ell\oti W_\ell 
  \\{}\\[-.8em]
  & & \bigoplus_{\ell}f_\ell\mortimes g_\ell\, \mapsto \bigoplus_{\ell}f_\ell\oti g_\ell
  \\{}\\[-.5em]
  {\rm and} \qquad
  & R{:}\quad \calC \To \CC\,,\quad&
  M\mapsto \bigoplus_{i\in\calI}(M\oti U_i\Vee)\times U_i \,,
  \\{}\\[-.8em]
  & & \,f\, \mapsto\, \bigoplus_{i\in\calI}(f\oti\id_{U_i\Vee})\mortimes\id_{U_i}  \,.
  \eear
  \ee
As shown in \cite{koRu2}, $R$ is a right adjoint for $T$, with the adjunction isomorphisms
  \be\label{adjunction}
  \phi_X^M{:}\quad \Hom_\C(T(X),M) \,\congTo\, \Hom_{\CC}(X,R(M))
  \ee
given by
  \be
  \phi^M_X:\ f = \bigoplus_{\ell}f_\ell \,\longmapsto\,
  \bigoplus_{\ell}\,\bigoplus_{i\in\I}\,\sum_\beta
  \Big[ ( [ f_\ell \cir (\id_{V_\ell} \oti g_\beta) ] \oti \id_{U_i^\vee} )
  \circ (\id_{V_\ell} \oti b_{U_i}) \Big] \mortimes \bar g_{\bar\beta} 
  \ee
for $X \eq \bigoplus_{\ell}V_\ell\,{\times}\, W_\ell$, where
$f_\ell \iN \Hom_\C(V_\ell\oti W_\ell,M)$ and
the $\beta$-summation is over a basis $\{g_\beta\}$ of $\Hom_\C(U_i,W_\ell)$
and the dual basis $\{\bar g_{\bar\beta}\}$ of $\Hom_\C(W_\ell,U_i)$. 
Note that again this is independent of the choice of a basis.
\\[3pt]
Now observe that $T(X)\eq\unit\boti X$ and $R(\unit)\eq A$. Thus for
$M\eq\unit$ the prescription \erf{adjunction} furnishes an isomorphism
  \be
  \phi_X^\unit{:}\quad \Hom_\C(\unit\boti X,\unit) \,\congTo\, \Hom_{\CC}(X,A) \,.
  \ee
This shows that $A$ possesses the defining property of $\intEnd(\unit)$.
Henceforth we abbreviate $\phi_X\,{\equiv}\, \phi_X^\unit$.

\medskip\noindent
(2)~ By \cite{ostr}, the evaluation morphism $\ev\iN\Hom_\C(\unit\boti A,\unit)$ 
is obtained as $\ev\eq\phi_{\!A}^{\,-1}(\id_A)$. Now we have
  \be
  \phi_A (\bigoplus_{j\in\calI}d_{U_j})=\bigoplus_{i,j\in\calI}
  \sum_\beta\Big[ ([d_{U_j}\cir (\id_{U_j\Vee}\oti g_\beta )]\oti \id_{U_i\Vee})\cir
  (\id_{U_j\Vee}\oti b_{U_i})\Big] \mortimes \bar g_\beta
  \ee
with $\{g_\beta\}$ a basis of $\Hom_\C(U_i,U_j)$. Since 
$\Hom_\C(U_i,U_j)\,{\cong}\,\delta_{i,j}\,\complex$, only $i\eq j$ contributes
to the $i$-summation. Choosing $g_\beta\eq\id_{U_i}$ we get
  \be
  \phi_A (\bigoplus_{j\in\calI}d_{U_j}) 
  =\bigoplus_{i\in\calI}\Big[ (d_{U_i}\cir \oti\id_{U_i\Vee})\cir (\id_{U_i\Vee}\oti b_{U_i})\Big]\mortimes\id_{U_i}
  =\bigoplus_{i\in\calI}\id_{U_i\Vee}\mortimes\id_{U_i}=\id_A \,.
  \ee
Thus $\ev \eq \phi_{\!A}^{\,-1}(\id_A) \eq \bigoplus_{j\in\calI}d_{U_j}$.

\medskip\noindent
(3)~ The multiplication on $\intEnd^{(n)}(\unit)$ is obtained as the image of
$\ev\cir(\ev\boti\id_A)\cir\psi^{(n)}_{\unit,A,A}$ under $\phi_{\!A\otimes A}$. We compute
  \be
  \bearl
  \ev\cir(\ev\boti\id_A)\cir\psi^{(n)}_{\unit,A,A}\\ 
  \qquad
  =\bigoplus_{i,j\in\calI}\,(d_{U_i}\oti d_{U_j})\cir (D^{-n}_{U_i\Vee U_i,U_j\Vee}\oti\id_{U_j})\cir\\
  \qquad\quad
  (\id_{U_i\Vee}\oti c_{U_j\Vee,U_i}\oti\id_{U_j})\cir(D^n_{U_i\Vee,U_j\Vee}\oti\id_{U_iU_j})
  \\
  \qquad
  =\bigoplus_{i,j\in\calI}d_{U_j}\cir(D^{-n}_{\unit,U_j\Vee}\oti\id_{U_j})\cir (d_{U_i}\oti\id_{U_j\Vee U_j})
 \cir\\
\qquad\quad
 (\id_{U_i\Vee}\oti c_{U_j\Vee,U_i}\oti\id_{U_j})\cir (D^n_{U_i\Vee,U_j\Vee}\oti\id_{U_iU_j})
  \\
  \qquad
  =\bigoplus_{i,j\in\calI}\,(d_{U_i}\oti d_{U_j})\cir (\id_{U_i\Vee}\oti c_{U_j\Vee,U_i}\oti\id_{U_j})\cir
  (D^n_{U_i\Vee,U_j\Vee}\oti\id_{U_iU_j})
  \eear
  \ee
Here the second equality uses naturality of $D_{U,V}$, and the third one that 
$D_{\unit,U}\eq\id_U$. We then find that, 
in terms of the dual bases $\{f_\alpha\}$ and $\{\bar f_\alpha\}$ in \erf{f-fbar-basis},
  \be
  \bearll
  m^{(n)} \!\!\!& = \phi_{A\otimes A}(\ev\cir(\ev\boti\id_A)\cir\psi^{(n)}_{\unit,A,A})
  \\[.4em]
  &= \bigoplus_{i,j,k\in\calI}\sum_\alpha\Big[\Big(\big[(d_{U_i}\oti d_{U_j})\cir (\id_{U_i\Vee}\oti 
  c_{U_j\Vee,U_i}\oti\id_{U_j})\cir(D^n_{U_i\Vee,U_j\Vee}\oti\id_{U_i U_j}) ~
  \\[.4em] \multicolumn2r{ 
  \circ\, (\id_{U_i\Vee U_j\Vee}\oti \bar f_\alpha)\big]\oti\id_{U_k\Vee}\Big)
  \cir(\id_{U_i\Vee U_j\Vee}\oti b_{U_k})\Big]\mortimes f_\alpha \,. }
  \eear
  \ee
This is indeed the morphism given in \erf{multiplication-morphism}.
\end{proof}

Since the module categories $(\C,\boxtimes,\psi^{(n)})$ are mutually equivalent, all the algebras
$(A,m^{(n)},\eta)$ are Morita equivalent. But in fact we even have

\begin{lemma}\label{lem_all_iso}
All algebra structures $(A,m^{(n)},\eta)$ with $n\iN\zet$ are isomorphic.
\end {lemma}

\begin{proof}
The morphism
  \be
  \sigma := \bigoplus_{i\in\I}\, (\theta_{U_i\Vee} \mortimes \id_{U_i})
  \label{def_sigma}\ee
is an automorphism of $A$ as an object. As an easy consequence of 
naturality and the defining property \erf{tensortwist} of the twist,
$\sigma$ indeed constitutes an algebra isomorphism 
from $(A,m^{(n)},\eta)$ to $(A,m^{(n+1)},\eta)$.
\end{proof}

It follows from the proof of Lemma \ref{lem_all_iso} that we can
express $m^{(n)}$ entirely through $m^{(0)}$ and the twist:
  \be
  m^{(n)} = \sigma^n \circ m^{(0)} \circ (\sigma^{-n} \oti \sigma^{-n}) 
  \ee
with $\sigma$ as given in \erf{def_sigma}.


\section{Frobenius algebras}\label{sec:frob}

That the object $A\eq \intEnd(\one)$ has the structure of an algebra is
guaranteed by the results of \cite{ostr}. But in fact $A$ carries more
structure, namely the one of a symmetric special Frobenius algebra. 
A {\em Frobenius\/} algebra in a monoidal category \X\ is an algebra in \X\ that
is also a coalgebra, with the coproduct being a morphism of $A$-bimodules.
For \X\ sovereign, a {\em symmetric\/} Frobenius algebra 
$B\eq(B,m,\eta,\Delta,\eps)$ in \X\ is a Frobenius algebra for which the two 
isomorphisms $((\eps\cir m)\oti\id_{B^\vee}) \cir (\id_B\oti b_B)$ and
$(\id_{B^\vee}\oti(\eps\cir m)) \cir (\tilde b_B\oti\id_B)$ from $B$ to its 
dual $B^\vee$ are equal, and a symmetric {\em special\/} Frobenius algebra 
is a symmetric Frobenius algebra for which $m\cir\Delta \eq \id_B$ and 
$\eps\cir\eta \eq \dim(B)\,\id_\one$ with $\dim(B)\,{\ne}\,0$.

The global dimension $\Dim(\C)$ of \C\ is the sum of the squares of dimensions 
of its simple objects,
  \be
  \Dim(\C) := \sum_{i \in \I} \dim(U_i)^2 .
  \ee
By \cite[Thm.\,2.3\,\&\,Cor\,2.10]{etno}, $\dim(U_i)$ is real 
and non-zero, and $\Dim(\C) \,{\ge}\, 1$.
  
\begin{proposition}\label{A-is-ssFA}
{\rm(1)}
$(A,m^{(n)},\eta,\Delta^{(n)},\eps)$ with $m^{(n)}$ and $\eta$ as defined in 
Theorem \ref{thm:End1=A} and with
  \be\label{comultiplication-morphism}
\raisebox{-58pt}{
  \bp(420,158)
\setlength{\unitlength}{1.2pt}
  \put(0,59)     {$\dsty \Delta^{(n)} ~:=~ \bigoplus_{i,j,k\in\I}
                   \frac{\dim(U_i)\,\dim(U_j)}{\Dim(\C)\,\dim(U_k)}\, \sum_\alpha $}
  \put(153,10){ \bp(2.5,0) \scalebox{.48}{\bild{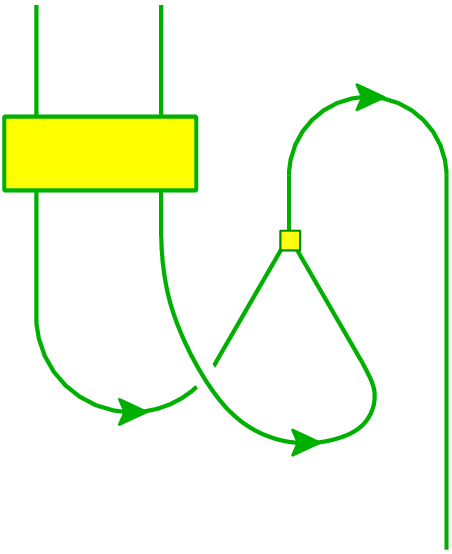}} \ep
   \put(-.5,111)  {\scriptsize $ U_i^\vee $}
   \put(24.6,112)  {\scriptsize $ U_{\!j}^\vee $}
   \put(77.5,-8)   {\scriptsize $ U_k^\vee $}
   \put(4.8,77.3)  {\scriptsize $ D_{U_i^\vee\!,U_{\!j}^\vee}^{-n}$}
   \put(56.4,58.1) {\scriptsize $ f_\alpha $}
  }
  \put(253,59)  {$ \mortimes $}
  \put(280,10){ \bp(0,0) \scalebox{.48}{\bild{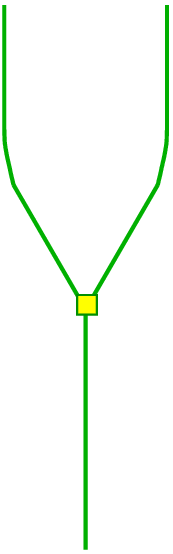}} \ep
   \put(-2.8,110)  {\scriptsize $ U_i $}
   \put(12.8,-8)   {\scriptsize $ U_k $}
   \put(28.2,110)  {\scriptsize $ U_j $}
   \put(19.1,45)   {\scriptsize $ \bar f_\alpha $}
  }
  \ep}
  \ee
(with the $\alpha$-summation as in \erf{multiplication-morphism}) and
  \be
  \eps := \Dim(\C)\, r_{\!A\succ \one\times\one}
  \ee
(with $r_{\!A\succ \one\times\one}$ the restriction morphism for $\,\one\,{\times}\,\one$ 
as a subobject of $A$) is a symmetric special Frobenius algebra.
\\[3mm]
{\rm(2)} The algebras $(A,m^{(n)},\eta,\Delta^{(n)},\eps)$, with
$n\iN\zet$, are all isomorphic as Frobenius algebras. 
\end{proposition}

\begin{proof}
(1)~\,For every $n$ consider the morphism
  \be
  \begin{split}
  \Phi^{(n)} &:= (d_A\oti\id_{A\Vee}) \cir (\id_A\oti m^{(n)}\oti\id_{A\Vee}) \cir
  (\tilde b_A\oti m^{(n)}\oti\id_{A\Vee}) \cir (\id_A\oti b_A) \\
&\in \Hom(A,A\Vee) \,.
\end{split}
  \ee
Starting from the definition of $m^{(n)}$ and using $D_{U,\one}\eq\id_U$ and 
$d_U\cir D_{U^\vee\!,U} \eq d_U\cir (\theta_{U^\vee}^2 \oti\id_{U})$,
it is straightforward to obtain the explicit expression
  \be
\raisebox{-45pt}{
  \bp(310,99)
  \put(0,44)       {$\dsty \Phi^{(n)} ~=~ \Dim(\C)\, \bigoplus_{i\in\I} \frac1{\dim(U_i)} $}
  \put(167,0){  {\bp(0,0) \scalebox{.3}{\bild{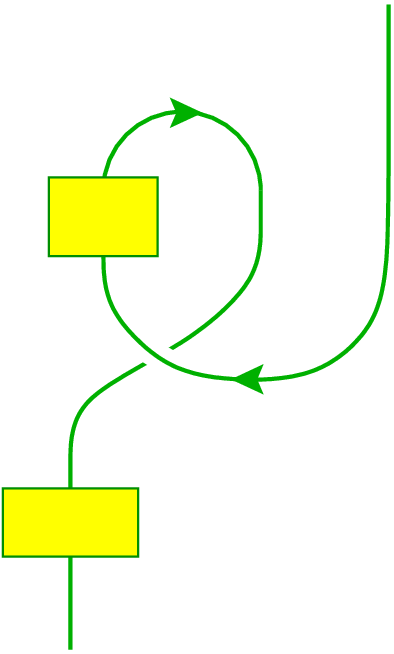}} \ep}
  \put(1.2,15.6)   {\scriptsize $ \theta^{-2n} $}
  \put(3.3,-8.8)   {\scriptsize $ U_i\Vee $}
  \put(7.7,61.1)   {\scriptsize $ \pi_i^{-1} $}
  \put(45.1,99)    {\scriptsize $ U_{\overline\imath}^{\,\vee\vee} $}
 }
  \put(240,44)     {$ \mortimes $}
  \put(266,0){  {\bp(0,0) \scalebox{.3}{\bild{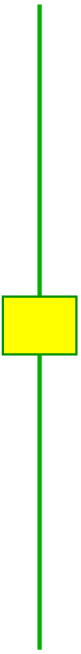}} \ep}
  \put(1.5,46.2)   {\scriptsize $ \pi_i $}
  \put(.1,99)      {\scriptsize $ U_{\overline\imath}^{\,\vee} $}
  \put(1.5,-8.8)   {\scriptsize $ U_i $}
 }
  \ep}
  \label{Phi}\ee
for these morphisms, where $\theta\,{\equiv}\,\theta_{U_i\Vee}$,
$\overline k\iN\I$ is the label for which $U_{\overline k} \,{\cong}\, U_k^\vee$, 
and $\pi_k\iN\Hom_\C(U_k^{},U_{\overline k}^\vee)$ is a choice of 
isomorphisms for each $k\in\I$ such that the compatibility condition (2.24) of 
\cite{fuRs4} between $\pi_k$ and $\pi_{\overline k}$ is satisfied. 
This shows in particular that $\Phi^{(n)}$ is invertible.
\\
Now according to Lemma 3.12 of \cite{fuRs4} invertibility of $\Phi^{(n)}$ implies
that $(A,m^{(n)},\eta)$ can be endowed with the structure of a symmetric special
Frobenius algebra, with the counit and coproduct given by $\eps^{(n)}\eq 
d_A\cir(\Phi^{(n)}\oti\eta) \eq d_A \cir (\id_A\oti m^{(n)}) \cir (\tilde b_A\oti \id_A)
$ and by 
  \be\bearll
  \Delta^{(n)} \!\! &= (m^{(n)}\oti{\Phi^{(n)}}^{-1}) \cir (\id_A\oti b_A)
  \\ 
  &= (d_A\oti{\Phi^{(n)}}^{-1}\oti{\Phi^{(n)}}^{-1})
  \cir (\id_A\oti m^{(n)}\oti\id_A\oti\id_A) \cir (\Phi^{(n)}\oti b_{A \otimes A}) \,,
  \eear \ee
respectively. Again by direct calculation one shows that these are precisely the
morphisms given in the proposition.
\\[.3em]
(2)~\,Two Frobenius algebras are isomorphic as Frobenius algebras if and only if
via one and the same isomorphism they are isomorphic both as algebras 
and as coalgebras. That all the algebras $(A,m^{(n)},\eta,\Delta^{(n)},\eps)$ are 
isomorphic as algebras, with an isomorphism given by the appropriate power of
$\sigma$ \erf{def_sigma}, has already been established in Lemma 
\ref{lem_all_iso}. Given the explicit form of the coproduct and counit, an entirely 
analogous argument shows that these morphisms constitute isomorphisms of coalgebras 
as well.
\end{proof}

Note that we can write $\Phi^{(n)} \eq \Phi^{(0)} \cir \sigma^{-2n}$ and
  \be
  \Delta^{(n)} = ( \sigma^n \oti \sigma^n) \circ \Delta^{(0)} \circ \sigma^{-n}
  \ee
with $\sigma$ as in \erf{def_sigma}. The coproduct $\Delta^{(0)}$ already
appeared in \cite[(2.57)]{koRu2}.


\section{Azumaya algebras}\label{sec:azumaya}

A {\em modular tensor category\/} \X\ is a premodular category for which
the symmetric $|\I|{\times}|\I|$-matrix $S$ with entries
\be
S_{ij}^{} \eq S_{00}\,(d_{X_j} \oti \tilde d_{X_i}) \cir [\, \id_{X_i^\vee} \oti
(c_{X_i,X_j}^{} \cir c_{X_j,X_i}^{}) \oti \id_{X_j^\vee} \,] \cir
(\tilde b_{X_j} \oti b_{X_i})
\ee
is invertible; here $\{X_i \,|\, i\iN\I \}$ 
denotes a set of representatives for the isomorphism classes of simple objects of \X, and 
$S_{00}\eq (\Dim(\C))^{-1/2} \,{>}\,0$.

If the category \C\ in Convention \ref{C-premodular} is even modular, then,
as we will see, the special symmetric Frobenius algebra $A$ from Proposition 
\ref{A-is-ssFA} has an additional property: it is Azumaya.

\begin{definition}
An {\em Azumaya category\/} \M\ over a braided monoidal category \X\
is a (right) module category \M\ over \X\ for which the two monoidal functors
$\alpha^\pm$ \erf{alphapm} from $(\X,\otimes^{\rm opp})$ to $\Endfun\X\M$ are equivalences.
\\[1pt]
An {\em Azumaya algebra\/} $B$ in a braided monoidal category \X\
is an algebra $B$ in \X\ such that the category of (left) $B$-modules is an
Azumaya category over \X.
\end{definition}

For $B$ a symmetric special Frobenius algebra in a modular tensor
category \X\ and $\X_B$ the category of left $B$-modules in \X, the endofunctor
category $\Endfun\X{\X_B}$ is monoidally equivalent to the category
\XBB\ of $B$-bimodules in \X. As a consequence, in this case the definition
of Azumaya algebra given here is equivalent to the one given in \cite{vazh}.

\medskip

For $B$ a symmetric special Frobenius algebra in a modular tensor category, 
consider the endomorphism
  \be
\raisebox{-42pt}{
  \bp(40,101)
  \put(0,0)   {\bp(0,0) \scalebox{.3}{\bild{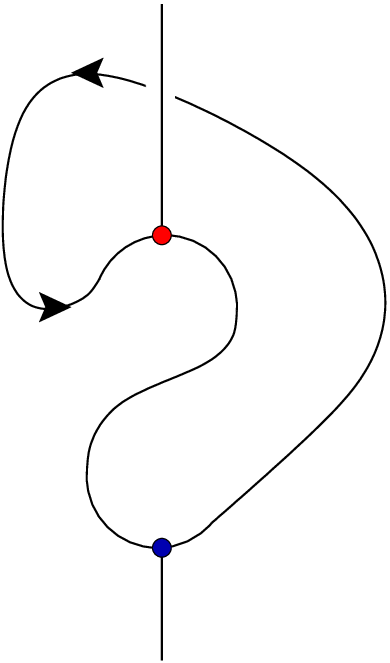}} \ep}
  \setlength{\unitlength}{.78pt}
  \put(-48.6,55.5) {$ P_{\!B}~:= $}
  \put(6.6,35)     {\scriptsize $ B $}
  \put(25.5,-8.8)  {\scriptsize $ B $}
  \put(26.1,125)   {\scriptsize $ B $}
  \put(73.2,67)    {\scriptsize $ B $}
  \ep}
  \label{P-def}\ee
~\\
which (e.g.\ by Lemma 5.2 of \cite{fuRs4}) is an idempotent.

\begin{lemma}\label{ssFA-is-Azu}
A symmetric special Frobenius algebra $B$ in a modular tensor category 
\X\ is Azumaya if and only if 
  \be
  P_{\!B} = \frac1{\dim(B)}\,\, \eta \circ \eps \,.
  \label{P-eta-eps}\ee
\end{lemma}

\begin{proof}
By \cite[Prop.\,2.30]{ffrs} the image of the idempotent $P_{\!B}$ of $B$ is the 
left center $C_l(B)$ as defined in \cite[Def.\,2.31]{ffrs}. Because of 
$\eps\cir\eta\eq\dim(B)\,\id_\one$, the equality \erf{P-eta-eps} is therefore 
equivalent to the statement that $C_l(B) \,{\cong}\, \one$.
Further, according to Proposition 2.36 of \cite{ffrs} one has
  \be
  \Hom_\X( C_l(B) \oti X, Y) \,\cong\, \Hom_{\Endfun\X\M}(\alpha^+_X,\alpha^+_Y) \,.
  \label{hom-Cl-iso} \ee
As we will now show, it follows that the braided induction 
functors $\alpha^\pm$ from \X\ to $\Endfun\X\M$ are equivalences iff $C_l(B)$ is trivial.
\\[.3em]
(i)~ Suppose that $\alpha^+$ is an equivalence of categories. Then by
\erf{hom-Cl-iso} we have $\Hom_\X( C_l(B) \oti X, Y)
     $\linebreak[0]$
{\cong}\,\Hom_\X(X,Y)$ 
for all $X,Y\iN\obj(\X)$, which implies that $C_l(B) \,{\cong}\, \one$.
\\[.3em]
(ii)\, Suppose that $C_l(B) \,{\cong}\, \one$. Then owing to the isomorphism
\erf{hom-Cl-iso} $\alpha^+$ is bijective on morphisms.
That $\alpha^+$ is also essentially surjective on objects, and hence
an equivalence, is seen as follows. Again by \erf{hom-Cl-iso},
$\alpha^+_{X_i}$ with $i\iN\I$ are mutually non-isomorphic simple objects of 
$\Endfun\X\M$. The same holds for $\{ \alpha^-_{X_i} \,|\, i\iN\I \}$. As a consequence,
each row of the $|\I|\,{\times}\,|\I|$-matrix with entries $ Z_{ij} \,{:=}\,
\dimc(\Hom_{\Endfun\X\M}(\alpha^+_{X_i},\alpha^-_{X_j})$ contains at most one non-zero
entry, and this non-zero entry, if present, must be equal to 1. It follows that 
$n \,{:=}\, \sum_{i,j \in \I} (Z_{ij})^2 \,{\le}\, |\I|$. On the other hand,
according to \cite[Sect.\,5.4\,\&\,Rem.\,5.19]{fuRs4}, $n$ is equal to the number
of isomorphism classes of simple objects in $\Endfun\X\M$. Since the $\alpha^+_{X_i}$ 
already provide $|\I|$ distinct simple objects, we must have $n \eq |\I|$.
Thus in particular $\alpha^+$ is essentially surjective, as claimed.
\\
The reasoning that $\alpha^-$ is an equivalence is analogous.
\end{proof}

\begin{proposition} \label{A-is-Azu}
If \C\ is a modular tensor category, then
the algebra $(A,m^{(n)},\eta)$ introduced in Theorem \ref{thm:End1=A} is Azumaya.
\end{proposition}

\begin{proof}
For $A$ as defined above, we have $\dim(A)\eq\Dim(\C)$ and \erf{P-eta-eps} reads
$P_{\!A} \eq e_{\one\times\one\prec A} \cir r_{\!A\succ \one\times\one}$,
which may be rewritten as $P_{\!A} \eq \bigoplus_{i\in \I} \delta_{i,0}^{}\, 
\id_{U_i^\vee\times U_i}$. Because of Lemma \ref{lem_all_iso} it suffices to 
establish this equality for the case $n\eq0$.

\medskip\noindent
Inserting the explicit form of the product and coproduct into \erf{P-def}, and rewriting the
resulting morphisms in the first factor of \CC\ with the help of the properties of the 
braiding (which amount to simple deformations of the graphical representation of the 
morphism), one shows that
  \be
\raisebox{-45pt}{
  \bp(370,105)
  \put(0,45)    {$\dsty P_{\!A}~=~\bigoplus_{i\iN\I}\sum_{j,k\in\I}
                 \frac{\dim(U_j)\,\dim(U_k)}{\Dim(\C)\,\dim(U_i)}\,\sum_{\alpha,\beta} $}
  \put(200,0)   {\bp(0,0) \scalebox{.3}{\bild{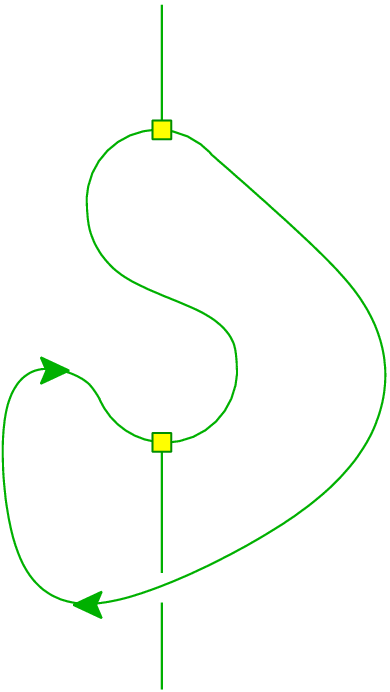}} \ep
    \put(25.8,82.8)  {\scriptsize $ \bar f_\alpha^\vee $}
    \put(19.8,42.7)  {\scriptsize $ {f'_\beta}^{\!\!\vee} $}
    \put(0.2,66)     {\scriptsize $ U_j^\vee $}
    \put(20.2,-8.8)  {\scriptsize $ U_i^\vee $}
    \put(18.8,103)   {\scriptsize $ U_i^\vee $}
    \put(45.2,67)    {\scriptsize $ U_k^\vee $}  }
  \put(265,45)  {$ \mortimes $}
  \put(290,0)   {\bp(0,0) \scalebox{.3}{\bild{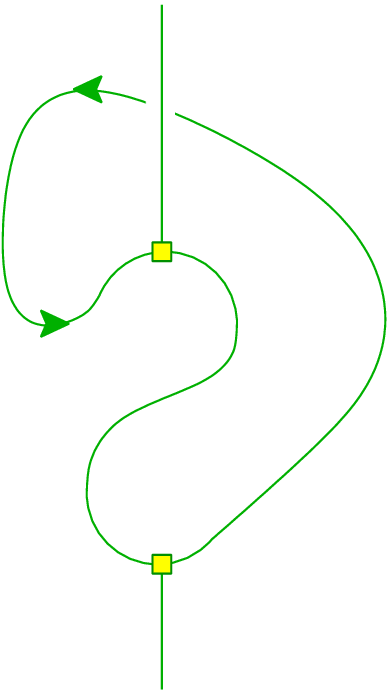}} \ep
    \put(19.4,55.1)  {\scriptsize $ f_\alpha $}
    \put(25.2,11.8)  {\scriptsize $ \bar f'_\beta $}
    \put(3,33)       {\scriptsize $ U_j $}
    \put(18.5,103)   {\scriptsize $ U_i $}
    \put(19.5,-8.8)  {\scriptsize $ U_i $}
    \put(48,72)      {\scriptsize $ U_k $}  }
  \ep}
  \ee
~\\
where the $\alpha$-summation is over a basis $\{f_\alpha\}$ of $\Hom_\C(U_k\oti U_j,U_i)$ 
and the corresponding dual basis of $\Hom_\C(U_i, U_k \oti U_j)$, 
and similarly for the $\beta$-summation. Next we note that the morphisms 
$u_{i;\alpha,\beta}^{(j,k)}\iN\End_\C(U_i\Vee)$ and $v_{i;\alpha,\beta}^{(j,k)}\iN\End_\C(U_i)$ 
that are represented by the two diagrams on the right 
hand side of this equality are multiples of $\id_{U_i\Vee}$ and of $\id_{U_i}$,
respectively. As a consequence, using properties of dual bases one can write
  \be
\raisebox{-45pt}{
  \bp(390,122)
  \put(0,55)    {$ \dsty u_{i;\alpha,\beta}^{(j,k)} \mortimes v_{i;\alpha,\beta}^{(j,k)} ~=~ 
                     \id_{U_i^\vee}^{} \mortimes [\,v_{i;\alpha,\beta}^{(j,k)} \cir w_{i;\alpha,\beta}^{(j,k)}
                     \,]\quad  $ with $\quad w_{i;\alpha,\beta}^{(j,k)}~:= $}
  \put(310,10)   {\bp(0,0) \scalebox{.3}{\bild{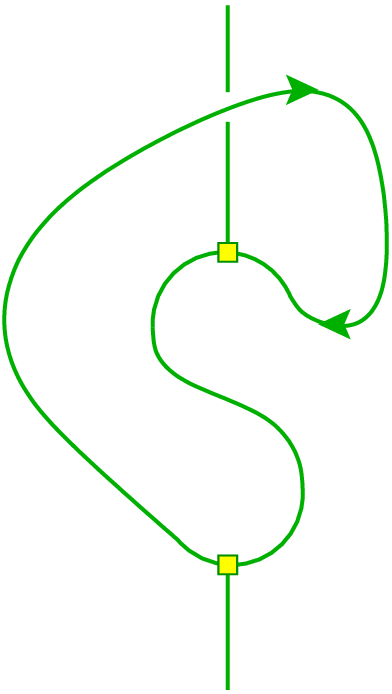}} \ep
    \put(23.1,11.1)  {\scriptsize $ \bar f_\alpha $}
    \put(22.7,68.3)  {\scriptsize $ f'_\beta $}
    \put(45.5,29)    {\scriptsize $ U_j $}
    \put(28.5,103)   {\scriptsize $ U_i $}
    \put(29.5,-8.8)  {\scriptsize $ U_i $}
    \put(-2.8,71)    {\scriptsize $ U_k $}  }
  \ep}
  \ee
Performing the $\beta$-summation leads to
  \be
\raisebox{-66pt}{
  \bp(360,157)
  \put(-60,76)    {$\dsty \sum_{\alpha,\beta}\,\id_{U_i^\vee}\mortimes
                    [v_{i;\alpha,\beta}^{(j,k)} \cir w_{i;\alpha,\beta}^{(j,k)}]
                    ~=~\sum_\alpha \id_{U_i^\vee}^{}\mortimes $}
  \put(160,10){
                {\bp(0,0) \scalebox{.3}{\bild{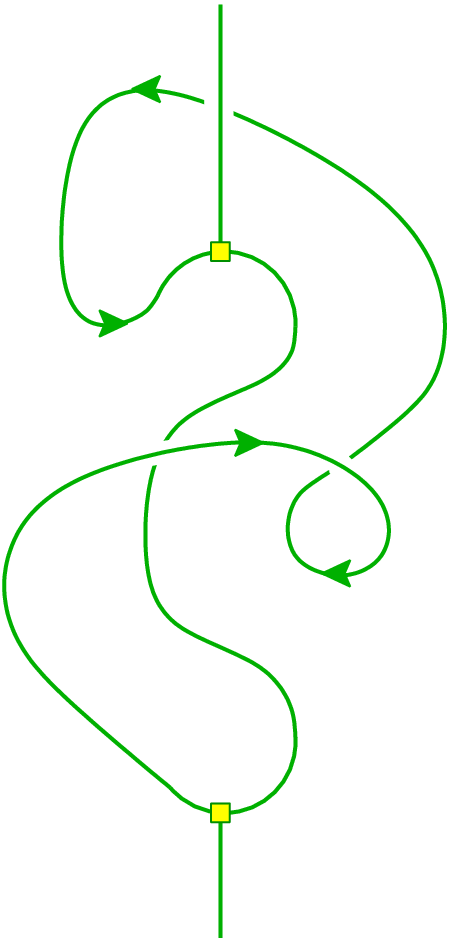}} \ep}
    \put(21.9,11.1)  {\scriptsize $ \bar f_\alpha $}
    \put(26.9,90.9)  {\scriptsize $ f_\alpha $}
    \put(28.7,-8.8)  {\scriptsize $ U_i $}
    \put(27.8,140)   {\scriptsize $ U_i $}
    \put(45.5,27)    {\scriptsize $ U_j $}
    \put(55.8,109)   {\scriptsize $ U_k $} }
  \put(242,76)  {$ =~\dsty\sum_\alpha \id_{U_i^\vee}^{}\mortimes $}
  \put(321,10){ {\bp(-8,0) \scalebox{.3}{\bild{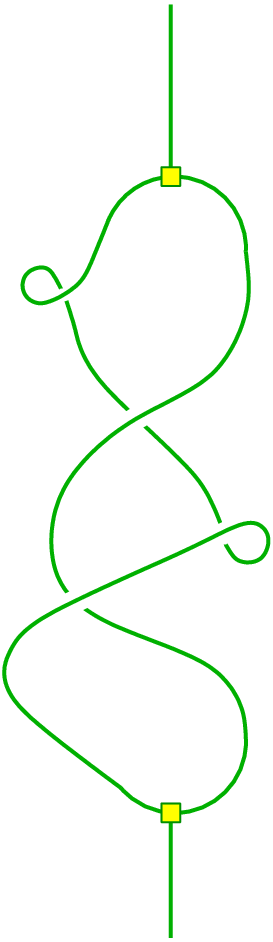}} \ep}
    \put(22.4,11.1)  {\scriptsize $ \bar f_\alpha $}
    \put(26.9,102)  {\scriptsize $ f_\alpha $}
    \put(29.5,-8.8)  {\scriptsize $ U_i $}
    \put(28.8,140)   {\scriptsize $ U_i $}
    \put(45.5,27)    {\scriptsize $ U_j $}
    \put(10.5,103)   {\scriptsize $ U_k $} }
  \ep}
  \label{qwerty}\ee
~\\
Using the naturality of the twist and the fact that the twist of a simple object is a multiple 
of the identity morphism, $\theta_{U_l} \eq \vartheta_l\,\id_{U_l}$, the second
factor on the right hand side of \erf{qwerty} can be rewritten as
  \be
  \begin{array}r \dsty
  f_\alpha\circ c_{U_j,U_k}^{} \cir c_{U_k,U_j}^{}\circ (\theta_{U_k}^2\oti \id_{U_j})
    \circ \bar f_\alpha 
  = \frac{\vartheta_k}{\vartheta_j}\, f_\alpha\circ c_{U_j,U_k}^{} \cir c_{U_k,U_j}^{} \circ
    (\theta_{U_k} \oti \theta_{U_j})\circ \bar f_\alpha ~\\{}\\[-.7em] \dsty
  = \frac{\vartheta_k}{\vartheta_j}\, f_\alpha \circ \theta_{U_k\otimes U_j}\circ \bar f_\alpha
  = \frac{\vartheta_k}{\vartheta_j}\, f_\alpha\cir\bar f_\alpha\circ \theta_{U_i}
  = \frac{\vartheta_k\,\vartheta_i}{\vartheta_j}\, \id_{U_i}\,.
  \eear
  \ee
After these manipulations the $\alpha$-summation has become trivial, yielding
just the dimension $\NN kji$ of the morphism space $\Hom_\C(U_k\oti U_j,U_i)$.
Summarising the calculations performed so far, we thus find that
  \be
  P_{\!A} = \bigoplus_{i\in \I} \xi_i\, \id_{U_i^\vee}^{} \mortimes \id_{U_i}^{}
  \qquad{\rm with}\qquad
  \xi_i = \sum_{j,k\in\I} \frac{\dim(U_j)\,\dim(U_k)}{\Dim(\C)\,\dim(U_i)}\,
  \frac{\vartheta_i\,\vartheta_k}{\vartheta_j}\, \NN kji \,.
  \label{P-xi}\ee
Now since \C\ is not only premodular, but even modular, then (see e.g.\ 
\cite[Ch.\,3.1]{BAki}) the dimensions $\NN ijk$ can be expressed through the 
$|\I|{\times}|\I|$-matrix $S$ 
as $\NN ijk \eq \sum_{l\in\I} S_{il}^{}\,S_{jl}^{}\,S_{lk}^*/S_{0l}^{}$,
and the matrices $S$ and $t$ with $t_{ij}^{} \eq \vartheta_i\,\delta_{i,j}^{}$
furnish a projective representation of the group SL$(2,\zet)$.
In particular we have $S\, t\, S \eq \gamma \, t^{-1} S\, t^{-1}\, C$ and 
$S\, t^{-1} S \eq \gamma^{-1} \, t\, S\, t$, where $\gamma \iN \complex^\times$ is 
a constant and $C$, the charge conjugation matrix, is the matrix with entries
$C_{k,l} \eq \delta_{k,\bar l}$. Together with $\dim(U_i)\eq S_{0i}/S_{00}$, this 
allows us to rewrite the numbers $\xi_i$ as
  \be
  \bearll \dsty
  \xi_i 
  \!\!&= \dsty S_{00}^{} \sum_{l\in\I} \frac{\vartheta_i\,S_{li}^*}{S_{0i}\,S_{0l}}\,
  (StS)_{0l}^{} \, (St^{-1}S)_{0l}^{}
  \\{}\\[-.7em] &\dsty
  = S_{00}^{} \sum_{l\in\I} \frac{\vartheta_i\,S_{li}^*}{S_{0i}\,S_{0l}}\,
  (t^{-1}St^{-1})_{0l}^{} \, (tSt)_{0l}^{}
  = S_{00}^{}\, \frac{\vartheta_i}{S_{0i}} \sum_{l\in\I} S_{il}^*\,S_{l0}^{}
  = \delta_{i,0}^{} \,,
  \eear
  \ee
so that \erf{P-xi} reduces to $P_{\!A} \eq \bigoplus_{i\in \I} \delta_{i,0}\, 
\id_{U_i^\vee}^{} \mortimes \id_{U_i}^{}$. Thus indeed $A$ is Azumaya.
\end{proof}

\begin{rem}
The product \erf{multiplication-morphism} and coproduct \erf{comultiplication-morphism} 
of $A$ do not make reference to the braiding in the second factor of \CbC. Denote by 
$\overline\C$ the modular tensor category obtained from \C\ by replacing the braiding 
$c$ and the twist $\theta$ by their inverses. Then $A$ is also a symmetric special
Frobenius algebra in $\C \,{\boxtimes}\, \overline\C$ (this is the situation considered in 
\cite{muge8,ffrs,koRu}). As an algebra in $\C \,{\boxtimes}\, \overline\C$, $A$ is commutative, 
i.e.\ $C_l(A) \eq A \eq C_r(A)$; in this sense it is `maximally non-Azumaya'.
\end{rem}


As an important consequence of the Azumaya property of $A$ we have

\begin{proposition}\label{prop-Zij,kl}
The spaces $\Hom_{\Endfun\CC\C}(\alpha^+_{U_i\times U_j},\alpha^-_{U_k\times U_l})$ have 
dimension
  \be
  \dimc(\Hom_{\Endfun\CC\C}(\alpha^+_{U_i\times U_j},\alpha^-_{U_k\times U_l}))
  = \delta_{i,l}\, \delta_{j,k} \,,
  \label{Zij,kl}\ee
for all $i,j,k,l \iN \I$.
\end{proposition}

\begin{proof}
Since \C\ is equivalent to the category of $A$-modules in \CC, the fact that $A$
is Azumaya implies \cite{fuRs11} that the $|\I|^2{\times}|\I|^2$-matrix $Z \eq (Z_{ij,kl})$ 
with $Z_{ij,kl}\eq \dimc(\Hom_{\Endfun\CC\C}(\alpha^+_{U_i\times U_j},\alpha^-_{U_k\times U_l}))$ 
is a permutation matrix, i.e.\ for any pair $(i,j)\iN\I{\times}\I$ there is precisely 
one pair $(k,l)$ such that $Z_{ij,kl}\eq 1$, while all other entries of $Z$ are zero.
Now according to Lemma \ref{lem_Gamma_M}, for any $i,j \iN \I$ the morphism space 
$\Hom_{\Endfun\CC\C}(\alpha^+_{U_i\times U_j},\alpha^-_{U_i\times U_j})$
is non-zero, i.e.\ $Z_{ij,ji}\eq 1$.
\end{proof}


\section{Permutation modular invariants}\label{sec:final}

We now address the application of our results to two-dimensional rational
conformal quantum field theory (CFT), to which we already alluded in the 
introduction.

A full local conformal field theory can be constructed from two ingredients:
a chiral conformal field theory, and one additional datum.
A chiral conformal field theory is, by definition, the following mathematical
object: a system of conformal blocks, i.e.\ sheaves over the moduli spaces of 
curves with marked points, with a Knizhnik-Zamolodchikov connection, 
constructed (see e.g.\ \cite{FRbe2,freN5}) from a conformal vertex
algebra. In the case of our interest, the vertex algebra \V\ is rational, which
amounts \cite{huan21} to the statement that its representation category 
\C\ is a modular tensor category. 

The additional datum needed for obtaining a full conformal field theory
is then a (Morita class of) symmetric special Frobenius algebra(s) $A$
in the representation category \C. The set of such Morita classes 
is in bijection to the set of consistent collections of correlation functions
defining a full local rational conformal field theory (see \cite{fjfrs2,koRu2}). 

To interpret the results of the preceding sections in this context,
suppose that the modular tensor category \C\ considered there is the
representation category of a rational vertex algebra \V. The monoidal category 
\CbC\ is then equivalent to the representation category of the vertex algebra 
$\V \,{\otimes_\complex}\, \V$. Further, the non-negative integers 
$Z_{ij,kl}(A) \eq \dimc(\Hom_{\Endfun\CC\C}(\alpha^+_{i,j},\alpha^-_{k,l}))$ 
are then the coefficients, in the standard basis of characters, of 
the {\em torus partition function\/} of the full CFT associated to the algebra
$A$ \cite{fuRs4}. It follows from Theorem 5.1 of \cite{fuRs4} that this
matrix obeys all requirements on a so-called {\em modular invariant\/}, as
formulated e.g.\ in Section V of \cite{gann16}. 

It should be appreciated, though, that not every quadratic matrix with non-negative
entries that obeys the axioms for a modular invariant describes the partition 
function of a full local conformal field theory. Counter examples even 
occur for the modular invariant induced by the duality on \C\ \cite{sosc}.
In view of this fact, the following application of the preceding
results is noteworthy. Take $A$ to be the symmetric special Frobenius algebra 
given in Proposition \ref{A-is-ssFA}. Then the result \erf{Zij,kl} amounts to the 
statement that the torus partition function of the full CFT determined by $A$
is given by the so-called $\zet_2$-{\em permutation modular invariant\/}. 
We thus arrive at

\begin{corollary}\label{cor_modinv}
The transposition modular invariant for a rational conformal field theory with
chiral symmetry $\V \,{\otimes_\complex}\, \V$ is physical, i.e.\ there exists
a consistent full local conformal field theory with this modular invariant 
as its torus partition function (and it can be
obtained via the construction of \cite{fuRs4}).
\end{corollary}

\begin{rem}
A word of warning is, however, in order: We do not have, at present, a 
classification of module categories over \CbC\ yielding the transposition modular 
invariant for the partition function, i.e.\ we cannot exclude that there exist 
symmetric special Frobenius algebras not Morita equivalent to the algebra $A$ 
considered here that give the same torus partition function as $A$.
\end{rem}

\begin{rem} \label{perm-vorb}
Proceeding along the lines of \cite{fuRs4} one can compute other correlation 
functions of the full CFT associated to the algebra $A$. One
of them is the partition function on an annulus. The coefficients of this
correlator in a standard basis are called the {\em annulus coefficients\/}.
As shown in \cite{fuRs4}, these are given by the integers
$\mathrm A_{X,M}^{~~~\,N} \,{:=}\, \dimc(\Hom_{A}(M\boti X,N))$ with 
$X$ a simple object and $M,N$ simple $A$-modules. For the algebra $A$ considered 
here, the induced $A$-modules $M_i \,{:=}\, \mathrm{Ind}_{A}(\one\,{\times}\,U_i)$ 
with $i\iN\I$ provide a complete set of pairwise nonisomorphic simple $A$-modules,
and a general induced $A$-module decomposes according to
$\mathrm{Ind}_{A}(U_i\,{\times}\,U_j) \,{\cong}\, \bigoplus_{k\in\I} \NN ijk\,M_k $
(compare Lemma 6.20 of \cite{ffrs}). As a consequence, the annulus coefficients are
  \be
  \mathrm A_{U_i\times U_j,M_k}^{~~~~~~~~\,M_l} = \NN{ij}k{\,l}
  \equiv \sum_{m\in\I} \NN ijm\,\NN mk{\,l}\equiv
  \dimc \Hom_\C(U_i\otimes U_j\otimes U_k,U_l) \,.
  \ee
Note that these non-negative integers are the dimensions of four-point conformal
blocks. The same dimensions appear in the fusion rules of permutation
orbifold theories \cite{bohs}. This does not come as a 
surprise: from the point of view of the permutation orbifold, the algebra
in Proposition \ref{A-is-ssFA} is Morita equivalent to the tensor unit,
and computing the annulus coefficients with the latter algebra,
one gets directly the fusion rules of the orbifold theory.
\end{rem}

At the same time, our results establish that the modular invariant
for the tensor product of an arbitrary finite number $N$ of identical
chiral conformal field theories obtained from {\em any\/} permutation 
$g\iN \mathfrak S_N$ is physical.
(Such modular invariants have been considered in string theory, see
e.g.\ \cite{fukS}.) This follows from the fact that every
permutation in $\mathfrak S_N$ is the product of transpositions, and that the
torus partition function obtained from a tensor product of special symmetric 
Frobenius algebras is given by the matrix product, see Proposition 5.3
of \cite{fuRs4}. We summarise these findings in

\begin{corollary}
Every permutation modular invariant for a rational conformal field theory with
chiral symmetry $\V \,{\otimes_\complex}\dots{\otimes_\complex}\, \V$ is physical,
i.e.\ there exists a consistent full local conformal field theory that has this 
modular invariant as its torus partition function 
(and it can be obtained via the construction of \cite{fuRs4}).
\end{corollary}

\noindent
For instance, it is easy to see that \C\ has the structure of a module
category over $\C{\boxtimes}\C{\boxtimes}\cdots{\boxtimes}\C$ ($m$ factors)
for which 
$$\intEnd(\unit)\,{\cong}\,\bigoplus_{i_1,i_2,...,i_m\in\I}
N_{i_1,i_2,...,i_m}\,U_{i_1}\,{\times}\,U_{i_2}\,{\times}\cdots{\times}\,
U_{i_m}$$ 
as an object of \C.

\smallskip

Let us also point out that already remark \ref{perm-vorb} hints at a relation 
with permutation orbifolds. Indeed, the permutation group $\mathfrak S_N$ acts
naturally on the $N$-fold tensor product $\V \,{\otimes_\complex}\dots
{\otimes_\complex}\, \V$ of vertex algebras. By general arguments one expects 
the invariant subalgebra, the so-called permutation orbifold chiral algebra, 
to be a rational vertex algebra as well. Its representation category, the 
modular tensor category $\C^{\rm orb}$ is of considerable interest, e.g.\ because 
its data enter all known approaches to the congruence subgroup conjecture
(see \cite{bant14,bant15}). This category cannot be obtained from $\C^{\boxtimes N}$
alone; rather, an $\mathfrak S_N$-equivariant modular category with $\C^{\boxtimes N}$ 
as its neutral component must be given. Our results provide the first steps in this direction:

\begin{rem}
Recall that, for \C\ a braided monoidal category and any $n,\hat n\iN\zet$,
\C\ together with $(\boxtimes,\psi^{(n)},r)$ as in Theorem \ref{thm_onlybraided}
and with $(\widehat\boxtimes,\widehat\psi^{(\hat n)},l)$ as in
Corollary \ref{cor_left} is both a left and a right module category over
\CtC. Actually, together with suitable natural mixed
associativity isomorphisms 
$$\widetilde\psi^{(n,\hat n)}_{X,M,Y}{:}~ 
X\botih(M\boti Y) \To (X\botih M)\boti Y$$
 it is even a {\em bimodule category\/}
over \CtC. The mixed associator is given by
  \be
  \widetilde\psi_{X,M,Y}^{(n,\hat n)} := \widehat\gamm^{\,\hat n}_{X,MY}
  \circ \big[ \id_X \botih \gamm^n_{M,Y} \big]
  \circ \big[ \widehat\gamm^{-\hat n}_{X,M} \boti \id_Y \big] \circ \gamm^{-n}_{XM,Y} 
  \ee
with
  \be
  \gamm_{M,U\times V}^{} := D_{M,U}^{\,-1} \oti \id_V \qquad{\rm and}\qquad
  \widehat\gamm_{U\times V,M}^{} := \id_U \oti D_{V,M}^{} \,.
  \ee
If \C\ has in addition a twist, then for all $n,\hat n\iN\zet$ the so obtained
structures of bimodule category on \C\ are equivalent.
This will be discussed in more detail elsewhere.
\end{rem}


\subsubsection*{Acknowledgements}
We thank Alexei Davydov for helpful discussions.
JF is partially supported by VR under project no.\ 621-2006-3343.
IR is partially supported by the EPSRC First Grant EP/E005047/1
and the Marie Curie network `Superstring Theory' (MRTN-CT-2004-512194). 
TB and CS are partially supported by the Collaborative Research Centre 676 ``Particles, 
Strings and the Early Universe - the Structure of Matter and Space-Time''.


 \newcommand\wb{\,\linebreak[0]} \def\wB {$\,$\wb}
 \newcommand\Bi[2]    {\bibitem[#2]{#1}}
 \newcommand\Erra[3]  {\,[{\em ibid.}\ {#1} ({#2}) {#3}, {\em Erratum}]}
 \newcommand\JO[6]    {{\em #6}, {#1} {#2} ({#3}), {#4--#5} }
 \newcommand\J[7]     {{\em #7}, {#1} {#2} ({#3}), {#4--#5} {{\tt [#6]}}}
 \newcommand\BOOK[4]  {{\em #1\/} ({#2}, {#3} {#4})}
 \newcommand\Prep[2]  {{\em #2}, preprint {#1} }

 \def\adma  {Adv.\wb in Math.}
 \def\anma  {Ann.\wb Math.}  
 \def\anop  {Ann.\wb Phys.}
 \def\apcs  {Applied\wB Cate\-go\-rical\wB Struc\-tures}  
 \def\atmp  {Adv.\wb Theor.\wb Math.\wb Phys.}
 \def\coma  {Con\-temp.\wb Math.}
 \def\comp  {Com\-mun.\wb Math.\wb Phys.}
 \def\fiic  {Fields\wb Inst.\wb Commun.}
 \def\ijmp  {Int.\wb J.\wb Mod.\wb Phys.\ A}
 \def\joac  {J.\wB Al\-ge\-bra\-ic\wB Com\-bin.}
 \def\joal  {J.\wB Al\-ge\-bra}
 \def\jpaa  {J.\wB Pure\wB Appl.\wb Alg.}
 \def\maan  {Math.\wb Annal.}
 \def\nupb  {Nucl.\wb Phys.\ B} 
 \def\pajm  {Pa\-cific\wB J.\wb Math.} 
 \def\phlb  {Phys.\wb Lett.\ B}
 \def\phrl  {Phys.\wb Rev.\wb Lett.}
 \def\pnas  {Proc.\wb Natl.\wb Acad.\wb Sci.\wb USA}
 \def\rvmp  {Rev.\wb Math.\wb Phys.}
 \def\sebo  {S\'emi\-nai\-re\wB Bour\-baki}
 \def\slnm  {Sprin\-ger\wB Lecture\wB Notes\wB in\wB Mathematics}
 \def\taac  {Theo\-ry\wB and\wB Appl.\wb Cat.} 
 \def\trgr  {Transform.\wb Groups}

\end{document}